\documentclass[authoryear, 11pt]{imsart}
\RequirePackage[OT1]{fontenc}
\RequirePackage{amsthm,amsmath,natbib}
\RequirePackage{hypernat}
\startlocaldefs
\numberwithin{equation}{section}
\theoremstyle{plain}

\endlocaldefs
\usepackage{psfrag, epsfig, graphicx}
\usepackage{amsfonts,amsmath,latexsym,amssymb}
\usepackage[active]{srcltx}
\usepackage{fullpage}

\newtheorem{lemma}{Lemma}
\newtheorem{theorem}{Theorem}
\newtheorem{proposition}{Proposition}
\newtheorem{assertion}{Assertion}

\newtheorem{assumption}{Assumption}
\newtheorem{remark}{Remark}
\newtheorem{definition}{Definition}
\newtheorem{corollary}{Corollary}

\newtheorem{assumption*}{Assumption}

\renewcommand{\kappa}{\varkappa}

\newcommand{\rd}{{\rm d}}
\newcommand{\e}{\varepsilon}

\newcommand{\cA}{{\cal A}}
\newcommand{\cB}{{\cal B}}
\newcommand{\cC}{{\cal C}}

\newcommand{\cF}{{\cal F}}

\newcommand{\cH}{{\cal H}}

\newcommand{\cK}{{\cal K}}
\newcommand{\cL}{{\cal L}}
\newcommand{\cM}{{\cal M}}

\newcommand{\cR}{{\cal R}}

\newcommand{\cU}{{\cal U}}

\newcommand{\cW}{{\cal W}}
\newcommand{\cX}{{\cal X}}

\newcommand{\BL}{\boldsymbol{L}}

\newcommand{\Bv}{\boldsymbol{v}}

\newcommand{\blg}{\boldsymbol{\mu}}

\newcommand{\bB}{\mathbb B}

\newcommand{\bE}{\mathbb E}
\newcommand{\bF}{\mathbb F}

\newcommand{\bH}{\mathbb H}

\newcommand{\bL}{{\mathbb L}}

\newcommand{\bN}{{\mathbb N}}

\newcommand{\bP}{{\mathbb P}}

\newcommand{\bR}{{\mathbb R}}

\newcommand{\bZ}{{\mathbb Z}}
\newcommand{\mA}{\mathfrak{A}}

\newcommand{\mP}{\mathfrak{P}}

\newcommand{\mM}{\mathfrak{M}}
\newcommand{\mU}{\mathfrak{U}}

\newcommand{\mJ}{\mathfrak{J}}
\newcommand{\mm}{\mathfrak{m}}

\newcommand{\mT}{\mathfrak{T}}

\newcommand{\mH}{\mathfrak{H}}

\newcommand{\mz}{\mathfrak{z}}
\newcommand{\ttt}{t}

\newcommand{\epr}{\hfill\hbox{\hskip 4pt
                \vrule width 5pt height 6pt depth 1.5pt}\vspace{0.5cm}\par}

\newcommand{\blh}{\boldsymbol{h}}

\newcommand{\bly}{\boldsymbol{y}}

\newcommand{\blL}{\boldsymbol{\Lambda}}

\begin{document}
\begin{frontmatter}
\title{Estimation in the convolution structure density model. Part I: oracle inequalities.}
\runtitle{Upper bounds}
\begin{aug}
\author[t1]
{\fnms{O.V.} \snm{Lepski}
\ead[label=e1]{oleg.lepski@univ-amu.fr}}
\author[t1]{\fnms{T.} \snm{Willer}
\ead[label=e2]{thomas.willer@univ-amu.fr}}
\thankstext{t1}{This work has been carried out in the framework of the Labex Archim\`ede (ANR-11-LABX-0033) and of the A*MIDEX project (ANR-11-IDEX-0001-02), funded by the "Investissements d'Avenir" French Government program managed by the French National Research Agency (ANR).}
\runauthor{O.V. Lepski and T. Willer}

\affiliation{Aix Marseille Univ, CNRS, Centrale Marseille, I2M, Marseille, France}

\address{Institut de Math\'ematique de Marseille\\
Aix-Marseille  Universit\'e   \\
 39, rue F. Joliot-Curie \\
13453 Marseille, France\\
\printead{e1}\\
\printead{e2}}
\end{aug}
\begin{abstract}

We study the problem of nonparametric estimation under $\bL_p$-loss, $p\in [1,\infty)$, in  the framework of the
convolution structure density model on $\bR^d$. This observation scheme is a generalization of two classical
statistical models, namely  density estimation under direct and indirect observations. In Part I the original
pointwise selection rule from a  family of "kernel-type" estimators is proposed.  For the selected estimator, we prove an $\bL_p$-norm oracle inequality and several of its consequences.
In Part II the problem of adaptive minimax estimation
under $\bL_p$--loss over the scale of anisotropic Nikol'skii classes is addressed.   We fully characterize the behavior of the minimax risk for different
relationships between regularity parameters and norm indexes in the definitions of the functional class and of the risk. We prove that the selection
rule proposed in Part I  leads to the construction of an optimally or nearly optimally (up to logarithmic factor)
adaptive estimator.

\end{abstract}
\begin{keyword}[class=AMS]
\kwd[]{62G05, 62G20}
\end{keyword}

\begin{keyword}
\kwd{deconvolution model}
\kwd{density estimation}
\kwd{oracle inequality}
\kwd{adaptive estimation}
\kwd{kernel estimators}
\kwd{$\bL_p$--risk}
\kwd{anisotropic Nikol'skii class}
\end{keyword}

\end{frontmatter}

\section{Introduction}

In the present paper we will investigate the  following observation scheme introduced in \cite{LW15}. Suppose that we observe i.i.d. vectors $Z_i\in\bR^d, i=1,\ldots, n,$ with a common probability density $\mathfrak{p}$ satisfying the following structural assumption
\vskip-0.4cm
\begin{equation}
\label{eq:convolution structure model}
\mathfrak{p}=(1-\alpha) f +\alpha [f\star g],\quad f\in\bF_g(R),\;\;\alpha\in[0,1],
\end{equation}
\vskip-0.1cm
\noindent where $\alpha\in[0,1]$ and  $g:\bR^d\to\bR$ are supposed to be known and $f:\bR^d\to\bR$ is the function to be estimated. We will call the observation scheme (\ref{eq:convolution structure model}) \textit{convolution structure density model}.

Here and later, for two functions $f,g\in\bL_1\big(\bR^d\big)$
\vskip-0.1cm
$$
\big[f\star g\big](x)=\int_{\bR^d}f(x-z)g(z)\nu_d(\rd z),\;\;x\in\bR^d,
$$
\vskip-0.1cm
\noindent and for any $\alpha\in[0,1]$, $g\in\bL_1\big(\bR^d\big)$ and $R> 1$,
\vskip-0.1cm
$$
\bF_g (R)=\Big\{f\in\bB_{1,d}(R): \;(1-\alpha) f +\alpha [f\star g]\in\mP\big(\bR^d\big)\Big\}.
$$
\vskip-0.1cm
\noindent Here $\mP\big(\bR^d\big)$ denotes the set of probability densities on $\bR^d$, $\bB_{s,d}(R)$ is the ball of radius $R>0$ in $\bL_s\big(\bR^d\big):=\bL_s\big(\bR^d,\nu_d\big), 1\leq s\leq\infty$ and  $\nu_d$ is the Lebesgue measure on $\bR^d$.

We remark that if one assumes additionally that $f,g\in\mP\big(\bR^d\big)$, this model can be interpreted as follows. The observations $Z_i\in\bR^d, i=1,\ldots, n,$ can be written as a sum of two independent random vectors, that is,
\begin{equation}
\label{eq:observation-scheme}
Z_i=X_i+\epsilon_iY_i,\quad i=1,\ldots,n,
\end{equation}
where $X_i, i=1,\ldots,n,$ are {\it i.i.d.} $d$-dimensional  random vectors with a common density $f$, to be estimated.
The noise variables  $Y_i, i=1,\ldots,n,$
are {\it i.i.d.} $d$-dimensional random vectors with a known common density $g$. At last $\e_i\in\{0,1\}, i=1,\ldots,n,$
are {\it i.i.d.} Bernoulli random variables with $\bP(\e_1=1)=\alpha$, where $\alpha\in [0,1]$ is supposed to be known.
The sequences  $\{X_i, i=1,\ldots,n\}$, $\{Y_i, i=1,\ldots,n\}$ and $\{\epsilon_i, i=1,\ldots,n\}$ are supposed to be mutually independent.

The observation scheme (\ref{eq:observation-scheme}) can be viewed as the generalization of two classical statistical models. Indeed, the case $\alpha=1$ corresponds to the standard deconvolution model $Z_i=X_i+Y_i,\; i=1,\ldots,n$. Another "extreme" case $\alpha=0$ corresponds to the direct observation scheme $Z_i=X_i,\; i=1,\ldots,n$.
The "intermediate" case $\alpha\in (0,1)$, considered for the first time in \cite{hesse}, can be treated as the mathematical modeling  of the following situation. One part of the data, namely  $(1-\alpha) n$, is observed without noise, while the other part is contaminated by additional noise. If the indexes corresponding to that first part were known, the density $f$ could be estimated using only this part of the data, with the accuracy corresponding to the direct case. The question we address now is: can one obtain the same accuracy if the latter information is not available? We will see that the answer to the aforementioned question is positive, but the construction of optimal estimation procedures is based upon ideas corresponding to the "pure" deconvolution model.

The convolution structure density model (\ref{eq:convolution structure model}) will be studied for an arbitrary $g\in\bL_1\big(\bR^d\big)$ and $f\in\bF_g(R)$.
Then, except in the case $\alpha=0$, the function $f$ is not necessarily a probability density.

We want to estimate $f$ using the observations $Z^{(n)}=(Z_1,\ldots,Z_n)$. By estimator, we mean any $Z^{(n)}$-measurable map $\hat{f}:\bR^n\to \bL_p\big(\bR^d\big)$. The accuracy of an estimator $\hat{f}$
is measured by the $\bL_p$--risk
\vskip-0.1cm
$$
 \cR^{(p)}_n[\hat{f}, f]:=\Big(\bE_f \|\hat{f}-f\|_p^p\Big)^{1/p},\;p\in [1,\infty),
$$
\vskip-0.1cm
\noindent where $\bE_f$ denotes the expectation with respect to the probability measure
$\bP_f$ of the observations $Z^{(n)}=(Z_1,\ldots,Z_n)$.
Also, $\|\cdot\|_p$, $p\in [1,\infty)$, is the $\bL_p$-norm on $\bR^d$ and without further mentioning we will assume that $f\in\bL_p\big(\bR^d\big)$.
The objective is to
construct an estimator of $f$ with a small $\bL_p$--risk.

\subsection{Oracle approach via local selection. Objectives of Part I}

Let $\cF=\big\{\hat{f}_\mathfrak{t}, \mathfrak{t}\in\mT\big\}$ be a family of estimators built from the observation $Z^{(n)}$. The goal is to propose a data-driven (based on $Z^{(n)}$) selection procedure from the collection $\cF$ and to establish for it an $\bL_p$-norm oracle inequality. More precisely, we want to construct a $Z^{(n)}$-measurable random map $\hat{\mathfrak{t}}:\bR^d\to \mT$
 and prove that for any $p\in [1,\infty)$ and $n\geq 1$
\begin{equation}
\label{eq:L_p-oracle-inequality}
\cR^{(p)}_n\big[\hat{f}_{\hat{\mathfrak{t}}(\cdot)}; f\big]\leq C_1\Big\|\inf_{\mathfrak{t}\in\mT}A_{n}\left(f,\mathfrak{t},\cdot\right)\Big\|_p+C_2n^{-\frac{1}{2}},\quad\forall f\in\bL_p\big(\bR^d\big).
\end{equation}
Here $C_1$ and $C_2$ are numerical constants which may depend on $d,p$ and $\mT$ only.

We call (\ref{eq:L_p-oracle-inequality}) an $\bL_p$-{\it norm oracle inequality obtained by local selection}, and in Part I
we provide with an explicit expression of the functional
$A_{n}(\cdot,\cdot,x), x\in\bR^d$ in the case where $\cF=\cF\big(\cH^d\big)$ is the family of "kernel-type" estimators parameterized by a collection of multi-bandwidths $\cH^d$.
The selection from the latter family is done pointwisely, i.e. for any $x\in\bR^d$, which allows to take into account the "local structure" of the function to be estimated. The $\bL_p$-norm oracle inequality is then obtained by the integration of the pointwise risk of the proposed estimator, which is a kernel estimator with the bandwidth being a multivariate random function. This, in its turn, allows us to derive different minimax adaptive results presented in Part II of the paper. They are obtained thanks to an unique $\bL_p$-norm oracle inequality.

Our selection rule presented in Section \ref{sec:susubsec-Pointwise-selection-rule-deconv} can be viewed as a generalization and modification of some statistical procedures proposed in \cite{{lepski-kerk}} and \cite{GL14}. As we mentioned above, establishing (\ref{eq:L_p-oracle-inequality}) is the main objective of Part I. We will see however  that although  $A_{n}(\cdot,\cdot,x), x\in\bR^d$ will be presented explicitly, its computation in particular problems is not a simple task. The main difficulty here is mostly related to the fact that (\ref{eq:L_p-oracle-inequality}) is proved without any assumption (except for the model requirements) imposed on the underlying function $f$. It turns out that under some nonrestrictive assumptions imposed on $f$, the obtained bounds can be considerably simplified, see Section \ref{sec:subscec:some consequences}. Moreover these new inequalities allow to better understand the methodology for obtaining minimax adaptive results by the use of the oracle approach.

\subsection{Adaptive estimation. Objectives of Part II}

Let $\bF$ be  a given subset of $\bL_p\big(\bR^d\big)$. For any estimator $\tilde{f}_n$, define its {\it maximal risk}  by
$
\cR^{(p)}_n\big[\tilde{f}_n; \bF\big]=\sup_{f\in\bF}\cR^{(p)}_n\big[\tilde{f}_n; f\big]
$
and its {\it minimax risk} on $\bF$ is given by
\vskip-0.6cm
\begin{equation}
\label{eq:minmax-risk}
\phi_n(\bF):=\inf_{\tilde{f}_n}\cR^{(p)}_n\big[\tilde{f}_n; \bF\big].
\end{equation}
\vskip-0.1cm
Here, the infimum is taken over all possible estimators. An estimator whose maximal risk is bounded, up to some constant factor, by $\phi_n(\bF)$, is called minimax on $\bF$.

Let $\big\{\bF_\vartheta,\vartheta\in\Theta\big\}$ be a collection of subsets of $\bL_p\big(\bR^d,\nu_d\big)$, where $\vartheta$ is a nuisance parameter which may have a very complicated structure.

The problem of adaptive estimation can be formulated as follows:
{\it is it possible to construct a single estimator $\hat{f}_n$
 which would be  simultaneously minimax on each class
 $\bF_\vartheta,\;\vartheta\in\Theta$, i.e.}
$$
 \limsup_{n\to \infty}\phi^{-1}_n(\bF_\vartheta)\cR^{(p)}_n\big[\hat{f}_n; \bF_\vartheta\big]<\infty,\;\;\forall \vartheta\in\Theta?
$$
We refer to this question as {\it the  problem of minimax adaptive
estimation over  the scale of }
$\{\bF_\vartheta,\;\vartheta\in\Theta \}$.
If such an estimator exists, we will call it optimally adaptive.


\smallskip

\noindent\textsf{From oracle approach to adaptation.} Let the oracle inequality (\ref{eq:L_p-oracle-inequality}) be established. Define
\vskip-0.1cm
$$
R_n\big(\bF_\vartheta\big)=\sup_{f\in\bF_\vartheta}\Big\|\inf_{\mathfrak{t}\in\mT}A_{n}\left(f,\mathfrak{t},\cdot\right)\Big\|_p+n^{-\frac{1}{2}},\quad \vartheta\in\Theta.
$$
\vskip-0.1cm
\noindent We immediately deduce from (\ref{eq:L_p-oracle-inequality}) that for any $\vartheta\in\Theta$
\vskip-0.1cm
$$
\limsup_{n\to\infty}R^{-1}_n\big(\bF_\vartheta\big)\cR^{(p)}_n\big[\hat{f}_{\hat{\mathfrak{t}}(\cdot)}; \bF_\vartheta\big]<\infty.
$$
\vskip-0.1cm
\noindent Hence, the minimax adaptive optimality of the estimator $\hat{f}_{\hat{\mathfrak{t}}(\cdot)}$ is reduced to the comparison of the normalization $R_n\big(\bF_\vartheta\big)$ with the minimax risk  $\phi_n(\bF_\vartheta)$.
Indeed,  if one  proves that for any $\vartheta\in\Theta$
\vskip-0.1cm
$$
\liminf_{n\to\infty}R_n\big(\bF_\vartheta\big)\phi^{-1}_n(\bF_\vartheta)<\infty,
$$
\vskip-0.1cm
\noindent then the estimator $\hat{f}_{\hat{\mathfrak{t}}(\cdot)}$ is {\it optimally adaptive} over the scale $\big\{\bF_\vartheta,\vartheta\in\Theta\big\}$.

\smallskip

\noindent\textsf{Objectives.}
In the framework of the {\it convolution structure density model}, we will be interested in adaptive estimation over the scale
\vskip-0.1cm
$$
\bF_\vartheta=\bN_{\vec{r},d}\big(\vec{\beta},\vec{L}\big)\cap\bF_g(R),\;\;\vartheta=\big(\vec{\beta},\vec{r},\vec{L},R\big),
$$
\vskip-0.1cm
where $\bN_{\vec{r},d}\big(\vec{\beta},\vec{L}\big)$ is the anisotropic Nikolskii
class (its exact definition will be presented in Part II). Here we only mention that for any $f\in\bN_{\vec{r},d}\big(\vec{\beta},\vec{L}\big)$ the coordinate $\beta_i$ of the vector  $\vec{\beta}=(\beta_1,\ldots,\beta_d)\in (0,\infty)^d$ represents the smoothness of $f$ in the direction $i$ and  the coordinate $r_i$ of the vector
$\vec{r}=(r_1,\ldots,r_d)\in [1,\infty]^d$ represents the index of the norm in which $\beta_i$ is measured. Moreover, $\bN_{\vec{r},d}\big(\vec{\beta},\vec{L}\big)$  is the intersection of balls in some semi-metric space and the vector $\vec{L}\in (0,\infty)^d$ represents the radii of these balls.

The aforementioned dependence on the direction is usually referred to {\it anisotropy} of the underlying function and the corresponding functional class.
The use of the integral norm in the definition of  the smoothness   is referred to {\it inhomogeneity} of the underlying function. The latter means that the function $f$ can be sufficiently smooth on some part of the  observation domain  and rather irregular on another part. Thus,  the adaptive  estimation over the scale $\big\{\bN_{\vec{r},d}\big(\vec{\beta},\vec{L}\big),\;\big(\vec{\beta},\vec{r},\vec{L}\big)\in(0,\infty)^d\times[1,\infty]^d\times(0,\infty)^d\big\}$
can be viewed as the adaptation to anisotropy  and inhomogeneity of the function to be estimated.

Additionally, we will consider
 $\bF_\vartheta=\bN_{\vec{r},d}\big(\vec{\beta},\vec{L}\big)\cap\bF_g(R)\cap\bB_{\infty,d}(Q),\; \vartheta=\big(\vec{\beta},\vec{r},\vec{L},R,Q\big)$.
It will allow us  to understand how the boundedness
of the underlying function may affect the accuracy of estimation.

The minimax adaptive estimation is a very active area of mathematical statistics, and the theory of adaptation was developed considerably over the past three decades. Several estimation procedures were proposed in various statistical models, such that
\textsf{Efroimovich-Pinsker method}, \cite{Pinsker, Efr}, \textsf{Lepski method},  \cite{lepski91} and \textsf{its generalizations}, \cite{{lepski-kerk}}, \cite{GL09},
\textsf{unbiased risk minimization}, \cite{golubev92}, \textsf{wavelet thresholding}, \cite{Donoho}, \textsf{model selection}, \cite{Barron-Birge, B-M}, \textsf{blockwise Stein method}, \cite{C99},
\textsf{aggregation of estimators}, \cite{nemirovski00}, \cite{Weg}, \cite{Tsyb03},  \cite{Gol}, \textsf{exponential weights}, \cite{Barron}, \cite{dal08},
\textsf{risk hull method}, \cite{cav-golubev}, among many others. The interested reader can find a very detailed overview as well as several open problems in adaptive estimation in the recent paper,
\cite{lepski15}.

As already said, the convolution structure density model includes itself the density estimation under direct and indirect observations. In Part II we compare in detail our minimax adaptive results to those already existing in both statistical models. Here we only mention that more developed results can be found in \cite{GL11}, \cite{GL14} (density model)  and in \cite{comte}, \cite{rebelles16} (density deconvolution).

\subsection{Assumption on the function $g$}

\noindent Later on for any $U\in\bL_1\big(\bR^d\big)$, let $\check{U}$ denote its Fourier transform, defined as $\check{U}(\ttt):=\int_{\bR^d}U(x)e^{-i\sum_{j=1}^d x_j t_j}\nu_d(dx), \ttt\in\bR^d.$
The selection rule from the family of kernel estimators, the $\bL_p$-norm oracle inequality as well as   the adaptive results presented in Part II
are established under the following condition.

\vskip-0.5cm
\begin{assumption}
\label{ass1:ass-on-noise-upper-bound}
(1) if $\alpha\neq 1$ then there exists  $\e>0$ such that
\vskip-0.2cm
$$
\big|1 - \alpha +\alpha\check{g}(t)\big|\geq \e, \quad \forall t\in\bR^d;
$$
\vskip-0.1cm
(2) if $\alpha=1$ then there exists  $\vec{\mu}=(\mu_1,\ldots,\mu_d)\in (0,\infty)^{d}$ and $\Upsilon_0>0$ such that
\vskip-0.3cm
$$|\check{g}(t)|\geq \Upsilon_0\prod_{j=1}^d(1+\ttt^2_j)^{-\frac{\mu_j}{2}},\quad\forall \ttt=(\ttt_1,\ldots,\ttt_d)\in\bR^d.
$$
\vskip-0.1cm
\end{assumption}

\vskip-0.1cm
Remind that the following assumption is well-known in the literature:
\vskip-0.1cm
$$
 \Upsilon_0\prod_{j=1}^d(1+t^2_j)^{-\frac{\mu_j}{2}}\leq |\check{g}(t)|\leq \Upsilon\prod_{j=1}^d(1+\ttt^2_j)^{-\frac{\mu_j}{2}},\quad\forall t\in\bR^d.
$$
\vskip-0.2cm
\noindent It is referred to as a {\it moderately ill-posed} statistical problem. In particular, the assumption is satisfied for the centered multivariate Laplace law.

Note that Assumption \ref{ass1:ass-on-noise-upper-bound} (1) is very weak and it is  verified for many distributions, including centered  multivariate Laplace and Gaussian ones.
Note also that this assumption always holds with $\e=1-2\alpha$ if $\alpha<1/2$.
Additionally, it holds with $\e=1-\alpha$ if $\check{g}$ is a real positive function. The latter is true, in particular, for any probability law obtained by an even number
of convolutions of a symmetric distribution with itself.

\section{Pointwise selection rule and $\bL_p$-norm oracle inequality}

To present our results in an unified way, let us define $\vec{\blg}(\alpha)=\vec{\mu}$, $\alpha=1$, $\vec{\blg}(\alpha)=(0,\ldots,0)$, $\alpha\in [0,1)$.
Let $K:\bR^d\to\bR$ be a continuous function belonging to $\bL_1\big(\bR^d\big)$, $\int_{\bR}K=1$, and such that its Fourier transform $\check{K}$ satisfies the following condition.

\begin{assumption}
\label{ass:on-kernel-deconvolution}
\vskip-0.2cm
There exist $\mathbf{k}_1>0$ and $\mathbf{k}_2>0$ such that
\vskip-0.5cm
\begin{eqnarray}
\int_{\bR^d}\big|\check{K}(t)\big|\prod_{j=1}^d(1+t^2_j)^{\frac{\blg_j(\alpha)}{2}}\rd t\leq \mathbf{k}_1,\quad
\int_{\bR^d}\big|\check{K}(t)\big|^2\prod_{j=1}^d(1+t^2_j)^{\blg_j(\alpha)}\rd t\leq \mathbf{k}^2_2.&&
\end{eqnarray}
\end{assumption}

\vskip-0.2cm
Set $\cH=\big\{e^{k},\; k\in\bZ\big\}$ and let
$
\cH^d=\big\{\vec{h}=(h_1,\ldots,h_d):\; h_j\in\cH, j=1,\ldots,d\big\}.
$
Define for any $\vec{h}=(h_1,\ldots,h_d)\in\cH^d$
\vskip-0.4cm
$$
K_{\vec{h}}(t)=V^{-1}_{\vec{h}}K\big(t_1/h_1,\ldots,t_d/h_d\big),\; t\in\bR^d,\quad V_{\vec{h}}=\prod_{j=1}^dh_j.
$$
\vskip-0.2cm
Later on for any $u,v\in\bR^d$ the operations and relations $u/v$, $uv$, $u\vee v$,$u\wedge v$,
$u\geq v$, $au, a\in\bR,$ are understood in coordinate-wise sense. In particular $u\geq v$ means that $u_j\geq v_j$ for any $j=1,\ldots, d$.

\subsection{Pointwise selection rule from the family of kernel estimators}
\label{sec:susubsec-Pointwise-selection-rule-deconv}

For any $\vec{h}\in (0,\infty)^d$ let  $M\big(\cdot,\vec{h}\big)$
satisfy the operator equation
\vskip-0.5cm
\begin{eqnarray}
\label{eq:def-kernel-M}
K_{\vec{h}}(y)=(1-\alpha)M\big(y,\vec{h}\big)+\alpha\int_{\bR^d}g(t-y)M\big(t,\vec{h}\big)\rd t,\quad  y\in\bR^d.&&
\end{eqnarray}
\vskip-0.3cm
\noindent For any $\vec{\mathrm{h}}\in\cH^d$ and $x\in\bR^d$ introduce the estimator
$
\widehat{f}_{\vec{\mathrm{h}}}(x)=n^{-1}\sum_{i=1}^n M\big(Z_i-x,\vec{\mathrm{h}}\big).
$


Our first goal is to propose for any given $x\in\bR^d$ a data-driven selection rule from the family of kernel estimators $\cF\big(\cH^d\big)=\big\{\widehat{f}_{\vec{\mathrm{h}}}(x),\;\vec{\mathrm{h}}\in\cH^d \big\}$.
Define for any $\vec{\mathrm{h}}\in\cH^d$
\vskip-0.5cm
\begin{gather*}
\widehat{U}_n\big(x,\vec{\mathrm{h}}\big)=\sqrt{\frac{2\lambda_n\big(\vec{\mathrm{h}}\big)\widehat{\sigma}^2\big(x,\vec{\mathrm{h}}\big)}{n}}+\frac{4 M_\infty\lambda_n\big(\vec{\mathrm{h}}\big)}{3n\prod_{j=1}^d\mathrm{h}_j(\mathrm{h}_j\wedge 1)^{\blg_j(\alpha)}},\quad \widehat{\sigma}^2\big(x,\vec{\mathrm{h}}\big)=\frac{1}{n}\sum_{i=1}^nM^2\big(Z_i-x,\vec{\mathrm{h}}\big);
\\
\lambda_n\big(\vec{\mathrm{h}}\big)=4\ln(M_\infty)+6\ln{(n)}+(8p+26)\sum_{j=1}^d\big[1+\blg_j(\alpha)\big]\big|\ln(\mathrm{h}_j)\big|;
\\
M_\infty=\big[(2\pi)^{-d}\big\{\e^{-1}\big\|\check{K}\big\|_1\mathrm{1}_{\alpha\neq 1}+\Upsilon_0^{-1}\mathbf{k}_1\mathrm{1}_{\alpha=1}\big\}\big]\vee 1.
\end{gather*}

\vskip-0.5cm
\paragraph{Pointwise selection rule} Let $\bH$ be an arbitrary subset of $\cH^d$. For any $\vec{h}\in\bH$ and $x\in\bR^d$ introduce
\vskip-0.5cm
\begin{eqnarray*}
&&\widehat{\cR}_{\vec{h}}(x)=\sup_{\vec{\eta}\in\bH}\Big[\big|\widehat{f}_{\vec{h}\vee\vec{\eta}}(x)-\widehat{f}_{\vec{\eta}}(x)\big|
-4\widehat{U}_n\big(x,\vec{h}\vee\vec{\eta}\big)-4\widehat{U}_n\big(x,\vec{\eta}\big)\Big]_+;
\\
&&\widehat{U}^*_n\big(x,\vec{h}\big)=\sup_{\vec{\eta}\in\bH:\; \vec{\eta}\geq\vec{h}}\widehat{U}_n\big(x,\vec{\eta}\big),
\end{eqnarray*}
\vskip-0.3cm
\noindent and define
\vskip-0.5cm
\begin{eqnarray}
\label{eq:pointwise-rule-deconvolution}
\vec{\mathbf{h}}(x)=\arg\inf_{\vec{h}\in\bH}\Big[\widehat{\cR}_{\vec{h}}(x)+8\widehat{U}^*_n\big(x,\vec{h}\big)\Big].
\end{eqnarray}
\vskip-0.3cm
\noindent Our final estimator is $\widehat{f}_{\vec{\mathbf{h}}(x)}(x),\;x\in\bR^d$ and we will call (\ref{eq:pointwise-rule-deconvolution}) the {\it pointwise selection rule}.

Note that the estimator $\widehat{f}_{\vec{\mathbf{h}}(\cdot)}(\cdot)$ does not necessarily belong to the collection $\big\{\widehat{f}_{\vec{\mathrm{h}}}(\cdot),\;\vec{\mathrm{h}}\in\cH^d \big\}$ since the multi-bandwidth $\vec{\mathbf{h}}(\cdot)$ is a $d$-variate function, which is not necessarily constant on $\bR^d$. The latter fact allows to take into account the "local structure" of the function to be estimated. Moreover, $\vec{\mathbf{h}}(\cdot)$ is chosen with respect to the observations, and therefore it is a random vector-function.

\subsection{$\bL_p$-norm oracle inequality}
\label{sec:subsec-l_p-norm-inequal-deconv}
Introduce for any $x\in\bR^d$ and $\vec{h}\in\cH^d$
\vskip-0.5cm
\begin{eqnarray*}
&&U^*_n\big(x,\vec{h}\big)=\sup_{\vec{\eta}\in\cH^d:\;\vec{\eta}\geq \vec{h}}U_n\big(x,\vec{\eta}\big),\qquad S_{\vec{h}}(x,f)=\int_{\bR^d}K_{\vec{h}}(t-x)f(t)\nu_d(\rd t);
\end{eqnarray*}
\vskip-0.3cm
\noindent where we have put
\vskip-0.5cm
\begin{gather*}
U_n\big(x,\vec{\eta}\big)=\sqrt{\frac{2\lambda_n\big(\vec{\eta}\big)\sigma^2\big(x,\vec{\eta}\big)}{n}}+\frac{4 M_\infty\lambda_n\big(\vec{\eta}\big)}{3n\prod_{j=1}^d\eta_j(\eta_j\wedge 1)^{\blg_j(\alpha)}},
\quad
\sigma^2\big(x,\vec{\eta}\big)=\int_{\bR^d}M^2\big(t-x,\vec{\eta}\big) \mathfrak{p}(t)\nu_d(\rd t).
\end{gather*}
\vskip-0.3cm
\noindent For any $\bH\subseteq\cH^d$, $\vec{h}\in\bH$ and  $x\in\bR^d$ introduce also
\begin{equation}
\label{eq:def-B_{vec{h}}}
B^*_{\vec{h}}(x,f)=\sup_{\vec{\eta}\in\bH}\big|S_{\vec{h}\vee \vec{\eta}}(x,f)-S_{ \vec{\eta}}(x,f)\big|,\qquad B_{\vec{h}}(x,f)=\big|S_{\vec{h}}(x,f)-f(x)\big|.
\end{equation}

\begin{theorem}
\label{th:L_p-norm-oracle-deconv}
\vskip-0.2cm
\noindent Let Assumptions \ref{ass1:ass-on-noise-upper-bound} and \ref{ass:on-kernel-deconvolution} be fulfilled. Then for any $\bH\subseteq\cH^d$, $n\geq 3$ and $p\in [1,\infty)$,
\vskip-0.2cm
$$
\forall f\in\bF_g(R),\quad \cR^{(p)}_n[\widehat{f}_{\vec{\mathbf{h}}(\cdot)}, f]\leq\Big\|\inf_{\vec{h}\in\bH}\Big\{2B^*_{\vec{h}}(\cdot,f)+B_{\vec{h}}(\cdot,f)+49U^*_n\big(\cdot,\vec{h}\big)\Big\}\Big\|_p+\mathbf{C}_pn^{-\frac{1}{2}}.
$$
\end{theorem}

\vskip-0.2cm
\noindent The explicit expression for the constant $\mathbf{C}_p$ can be found in the proof of the theorem.


Later on we will pay attention to a special choice for the collection of multi-bandwidths, namely
\vskip-0.2cm
$$
\cH^d_{\text{isotr}}:=\big\{\vec{h}\in\cH^d:\; \vec{h}=(h,\ldots,h), \;\; h\in\cH\big\}.
$$
\vskip-0.1cm
\noindent More precisely, in Part II, the selection from the corresponding family of kernel estimators will be used for the adaptive estimation over the collection of isotropic Nikolskii classes. Note also that if $\bH=\cH^d_{\text{isotr}}$ then obviously for any $ \vec{h}=(h,\ldots,h)\in \cH^d_{\text{isotr}}$
\vskip-0.2cm
$$
B^*_{\vec{h}}(\cdot,f)\leq 2\sup_{\vec{\eta}\in\cH^d_{\text{isotr}}:\; \eta\leq h}B_{\vec{\eta}}(\cdot,f)
$$
\noindent and we come to the following corollary of Theorem \ref{th:L_p-norm-oracle-deconv}.
\begin{corollary}
\label{cor:th-L_p-norm-oracle-deconv}
Let Assumptions \ref{ass1:ass-on-noise-upper-bound} and \ref{ass:on-kernel-deconvolution} be fulfilled. Then for any  $n\geq 3$ and $p\in [1,\infty)$
\vskip-0.4cm
$$
\cR^{(p)}_n[\widehat{f}_{\vec{\mathbf{h}}(\cdot)}, f]\leq\bigg\|\inf_{\vec{h}\in\cH^d_{\text{isotr}}}\bigg\{5\sup_{\vec{\eta}\in\cH^d_{\text{isotr}}\;:\; \eta\leq h}B_{\vec{\eta}}(\cdot,f)+49U^*_n\big(\cdot,\vec{h}\big)\bigg\}\bigg\|_p+\mathbf{C}_pn^{-\frac{1}{2}},
\quad\forall f\in\bF_g(R).
$$
\vskip-0.0cm
\end{corollary}
The oracle inequality proved in Theorem \ref{th:L_p-norm-oracle-deconv} is particularly useful since it does not require any assumption on the underlying function $f$ (except for the restrictions ensuring the existence of the model and of the risk). However, the quantity appearing in the right hand side of this inequality, namely
$$
\Big\|\inf_{\vec{h}\in\bH}\Big\{2B^*_{\vec{h}}(\cdot,f)+B_{\vec{h}}(\cdot,f)+49U^*_n\big(\cdot,\vec{h}\big)\Big\}\Big\|_p
$$
\noindent is not easy to analyze. In particular, in order to use the result of Theorem \ref{th:L_p-norm-oracle-deconv}   for adaptive estimation, one has to be able to compute
$$
\sup_{f\in\bF}\Big\|\inf_{\vec{h}\in\bH}\Big\{2B^*_{\vec{h}}(\cdot,f)+B_{\vec{h}}(\cdot,f)+49U^*_n\big(\cdot,\vec{h}\big)\Big\}\Big\|_p
$$
\noindent for a given class $\bF\subset \bL_p\big(\bR^d\big)\cap \bF_g(R)$ with either $\bH=\cH^d$ or $\bH=\cH^d_{\text{isotr}}$. It turns out that under some nonrestrictive  assumptions imposed on $f$, the obtained bounds can be considerably simplified. Moreover, the new inequality obtained below will allow us to better understand the way for proving adaptive results.

\subsection{Some consequences of Theorem \ref{th:L_p-norm-oracle-deconv}}
\label{sec:subscec:some consequences}

Thus, furthermore we will assume that $f\in\bF_{g,\mathbf{u}}(R,D)\cap\bB_{\mathbf{q},d}(D)$, $\mathbf{q},\mathbf{u}\in[1,\infty], D>0,$ where
$$
\bF_{g,\mathbf{u}} (R,D):=\Big\{f\in\bF_g(R): \;(1-\alpha) f +\alpha [f\star g]\in\bB^{(\infty)}_{\mathbf{u},d}(D)\Big\},
$$
\noindent and $\bB^{(\infty)}_{\mathbf{u},d}(D)$ denotes the ball of radius $D$ in the weak-type space  $\bL_{\mathbf{u},\infty}\big(\bR^d\big)$, i.e.
$$
\bB^{(\infty)}_{\mathbf{u},d}(D)=\big\{\lambda:\bR^d\to\bR: \; \|\lambda\|_{\mathbf{u},\infty}< D\big\},\; \|\lambda\|_{\mathbf{u},\infty}=\inf\big\{C:\; \nu_d\big(x: |T(x)|>\mz\big)\leq C^\mathbf{u}\mz^\mathbf{-u},\;\forall \mz>0\big\}.
$$
\noindent As usual $\bB^{(\infty)}_{\mathbf{\infty},d}(D)=\bB_{\infty,d}(D)$ and obviously $\bB^{(\infty)}_{\mathbf{u},d}(D)\supset\bB_{\mathbf{u},d}(D)$.
Note also that $\bF_{g,\mathbf{1}} (R,D)=\bF_{g} (R)$ for any $D\geq 1$.
It is worth noting that the assumption $f\in\bF_{g,\mathbf{u}} (R,D)$ simply means that the common density of the observations $\mathfrak{p}$ belongs to $\bB^{(\infty)}_{\mathbf{u},d}(D)$.

\begin{remark}
\label{rem:before-sec-inequality}
\vskip-0.2cm
It is easily seen that $\bF_{g,\infty} \big(R, R\|g\|_\infty\big)=\bF_g(R)$ if $\alpha=1$ and $\|g\|_\infty<\infty$. Note also that
$
\bF_{g,\mathbf{\infty}}(R, Q\|g\|_1)\supset \bF_g(R)\cap \bB_{\infty,d}(Q)
$
for any $\alpha\in[0,1]$ and $Q>0$.
\end{remark}

\subsubsection{\textsf{Oracle inequality over $\bF_{g,\mathbf{u}} (R,D)\cap\bB_{\mathbf{q},d}(D)$}}
For any $\vec{h}\in\cH^d$ and any $v>0$, let
$$
\cB_{\vec{h}}(\cdot,f)=2B^*_{\vec{h}}(\cdot,f)+B_{\vec{h}}(\cdot,f),\quad \cA(\vec{h},f,v)=\big\{x\in\bR^d:\; \cB_{\vec{h}}(x,f)\geq 2^{-1}v\big\},
$$
$$
F_n\big(\vec{h}\big)=\frac{\sqrt{\ln{n}+ \sum_{j=1}^d|\ln{h_j|}}}{\sqrt{n}\prod_{j=1}^dh^{\frac{1}{2}}_j(h_j\wedge 1)^{\blg_j(\alpha)}}, \quad G_n\big(\vec{h}\big)=\frac{\ln{n}+ \sum_{j=1}^d|\ln{h_j|}}{n\prod_{j=1}^dh_j(h_j\wedge 1)^{\blg_j(\alpha)}}.
$$
\noindent Furthermore let $\bH$ be either $\cH^d$ or $\cH^d_{\text{isotr}}$ and for any $v,z>0$ define
\begin{equation}
\label{eq:def-u(bly)}
\mH(v)=\big\{\vec{h}\in\bH:\; G_n\big(\vec{h}\big) \leq av\big\},\quad \mH(v,z)=\big\{\vec{h}\in\mH(v):\;
F_n\big(\vec{h}\big) \leq avz^{-1/2}\big\}.
\end{equation}
Here  $a>0$ is a numerical  constant  whose explicit expression is given in the beginning of Section \ref{sec:subsec-proof-concequences}.
Introduce for any $v>0$ and $f\in\bF_{g,\mathbf{u}} (R,D)$
\vskip-0.4cm
\begin{eqnarray*}
\Lambda(v,f)&=&\inf_{\vec{h}\in\mH(v)}\Big[\nu_d\big(\cA(\vec{h},f,v)\big)+v^{-2}F_n^2\big(\vec{h}\big)\Big];
\\
\Lambda(v,f,\mathbf{u})&=&\inf_{z\geq 2}\;\inf_{\vec{h}\in\mH(v,z)}\Big[\nu_d\big(\cA(\vec{h},f,v)\big)+z^{-\mathbf{u}}\Big];
\\
\Lambda_{p}(v,f,\mathbf{u})&=&\inf_{z\geq 2}\;\inf_{\vec{h}\in\mH(v,z)}\bigg[\int_{\cA(\vec{h},f,v)}\big|\cB_{\vec{h}
}(x,f)\big|^p\nu_d(\rd x)+v^{p}z^{-\mathbf{u}}\bigg].
\end{eqnarray*}

\begin{remark}
\label{rem:existence-of-Lambdas}
Note that $\mH(v)\neq\emptyset$ and $\mH(v,z)\neq\emptyset$ whatever the values of $v>0$ and $z\geq 2.$
Indeed, for any $v>0$ and $z>2$ one can find $b>1$ such that
$$
\big(\ln{n}+ d\ln{b}\big)(nb^d)^{-1}\leq \big[a^2v^2z^{-1}\big]\wedge av.
$$
The latter means that $\vec{b}=(b,\ldots,b)\in\mH(v,z)\cap\mH(v)$.
Thus, we conclude that the quantities $\Lambda(v,f)$, $\Lambda(v,f,\mathbf{u})$ and $\Lambda_{p}(v,f,\mathbf{u})$ are well-defined for all $v>0$.

\end{remark}

Also, It is easily seen that for any $v>0$ and $f\in\bF_{g,\mathbf{\infty}} (R,D)$
\begin{equation}
\label{eq:case-mathbf(u)=infty}
\Lambda(v,f,\infty)=\inf_{\vec{h}\in\mH(v,2)}\nu_d\big(\cA(\vec{h},f,v)\big),\quad \Lambda_{p}(v,f,\infty)=\inf_{\vec{h}\in\mH(v,2)}\int_{\cA(\vec{h},f,v)}\big|\cB_{\vec{h}
}(x,f)\big|^p\nu_d(\rd x).
\end{equation}
Put at last for any $v>0$, $l_\bH(v)=v^{p-1}(1+|\ln{(v)}|)^{t(\bH)}$, where $t(\bH)=d-1$ if $\bH=\cH^d$ and $t(\bH)=0$ if  $\bH=\cH^d_{\text{isotr}}$.

\begin{theorem}
\label{th:consec-th1-unbounded}
Let the assumptions of Theorem \ref{th:L_p-norm-oracle-deconv} be fulfilled and let $K$ be a compactly supported function.
Then for any $n\geq 3$, $p>1, \mathbf{q}>1, R>1, D>0,  0<\underline{\Bv}\leq\overline{\Bv}<\infty, \mathbf{u}\in(p/2,\infty], \mathbf{u}\geq \mathbf{q}$ and any $f\in\bF_{g,\mathbf{u}} (R,D)\cap\bB_{\mathbf{q},d}(D)$
\begin{equation*}
\cR^{(p)}_n[\widehat{f}_{\vec{\mathbf{h}}(\cdot)}, f]\leq C^{(1)}\bigg[l_\bH(\underline{\Bv})+\int_{\underline{\Bv}}^{\overline{\Bv}}v^{p-1}\{\Lambda(v,f)\wedge\Lambda(v,f,\mathbf{u})\}\rd v+
\Lambda_{p}(\overline{\Bv},f,\mathbf{u})
\bigg]^{\frac{1}{p}}+\mathbf{C}_pn^{-\frac{1}{2}}.
\end{equation*}
\end{theorem}
Here $C^{(1)}$  is a universal  constant independent of $f$ and $n$. Its explicit expression  can be found in the proof
of the  theorem. We remark also that only this constant depends on $\mathbf{q}$.

The result announced in Theorem \ref{th:consec-th1-unbounded} suggests a way for establishing  minimax and minimax adaptive properties of the pointwise selection rule given in (\ref{eq:pointwise-rule-deconvolution}). For a given $\bF\subset\bF_{g,\mathbf{u}} (R,D)\cap\bB_{\mathbf{q},d}(D)$ it mostly consists in finding a careful estimate for
$$
S\big(\vec{h},\mz\big):=\sup_{f\in\bF}\nu_d\big(x\in\bR^d:\; \cB_{\vec{h}}(x,f)\geq \mz\big),\quad \forall \vec{h}\in\bH, \;\; \forall\mz>0.
$$
The choice of $\underline{\Bv},\overline{\Bv}>0$ is a delicate problem and it  depends on $S(\cdot,\cdot)$.

In the next section we present several results concerning  some useful upper estimates for the quantities
$$
\sup_{f\in\bF}\{\Lambda(v,f)\wedge\Lambda(v,f,\mathbf{u})\},\quad \sup_{f\in\bF}\Lambda_{p}(v,f,\mathbf{u}),\; v>0.
$$
We would like to underline that these bounds will be established for an arbitrary $\bF$ and, therefore, they can be applied to the adaptation over different scales of functional classes. In particular, the results obtained below form the basis for our consideration in Part II.

\subsubsection{\textsf{Application to the minimax adaptive estimation}}

Our objective now is to bound from above
$
\sup_{f\in\bF}\cR^{(p)}_n[\widehat{f}_{\vec{\mathbf{h}}(\cdot)}, f]
$
for any $\bF\subset \bF_{g,\mathbf{u}} (R,D)\cap\bB_{\mathbf{q},d}(D)$.
All the results in this section will be proved under an additional condition imposed on the kernel $K$.

\begin{assumption}
\label{ass2:on-kernel-deconvolution-part1}
Let $\cK:\bR\to\bR$ be a compactly supported, bounded function and $\int \cK=1$. Then
\vskip-0.3cm
$$
K(x)=\prod_{j=1}^d \cK(x_j),\; \forall x\in\bR^d.
$$
\vskip-0.2cm
Without loss of generality we will assume that $\|\cK\|_\infty\geq 1$ and $\text{supp}(\cK)\subset[-c_\cK,c_\cK]$ with $c_\cK\geq 1$.
\end{assumption}

Introduce the following notations. Set for any $h\in\cH$,  $x\in\bR^d$ and $j=1,\ldots,d$
$$
b^*_{h,f,j}(x)=\left|\int_{\bR}\cK(u)f\big(x+uh\mathbf{e}_j\big)\nu_1(\rd u)-f(x) \right|,
\quad b_{h,f,j}(x)=\sup_{\eta\in\cH:\: \eta\leq h}b^*_{\eta,f,j}(x)
$$
where $(\mathbf{e}_1,\ldots,\mathbf{e}_d)$ denotes the canonical basis of $\bR^d$. For any $s\in[1,\infty]$ introduce
$$
\mathbf{B}^*_{j,s,\bF}(\mathbf{h})=
\sup_{f\in\bF}\sum_{h\in\cH:\: h\leq \mathbf{h}}\big\|b^*_{h,f,j}\big\|_s,
\quad \mathbf{B}_{j,s,\bF}(\mathbf{h})=\sup_{f\in\bF}\big\|b_{\mathbf{h},f,j}\big\|_{s},
\quad j=1,\ldots,d.
$$
Set for any $\vec{h}\in\cH^d$,
$v>0$  and  $j=1,\ldots,d$,
$$
J\big(\vec{h},v\big)=\big\{j\in\{1,\ldots,d\}:\; h_j\in\mathbf{V}_{j}(v)\big\},\;\;\;\mathbf{V}_{j}(v)=\big\{\mathbf{h}\in\cH:\;\; \mathbf{B}_{j,\infty,\bF}(\mathbf{h})\leq \mathbf{c}v\big\},
$$
where $\mathbf{c}=(20d)^{-1}\big[\max (2c_\cK\|\cK\|_\infty,\|\cK\|_1)\big]^{-d}$.
As usual the complement of $J\big(\vec{h},v\big)$ will be denoted by $\bar{J}\big(\vec{h},v\big)$. Furthermore, the summation over the empty set is supposed to be zero.

For any $\vec{s}=(s_1,\ldots,s_d)\in[1,\infty)^d$, $\mathbf{u}\geq 1$ and $v>0$ introduce
\begin{gather}
\label{eq4:def-blL}
\blL_{\vec{s}}(v,\bF,\mathbf{u})=\inf_{z\geq 2}\;\inf_{\vec{h}\in\mH(v,z)}\bigg[\sum_{j\in\bar{J}(\vec{h},v)} v^{-s_j}\big[\mathbf{B}_{j,s_j,\bF}\big(h_{j}\big)\big]^{s_j}+z^{-\mathbf{u}}\bigg];
\\
\label{eq1:def-blL}
\blL_{\vec{s}}\big(v,\bF\big)=\inf_{\vec{h}\in\mH(v)}\bigg[\sum_{j\in\bar{J}(\vec{h},v)} v^{-s_j}\big[\mathbf{B}_{j,s_j,\bF}\big(h_{j}\big)\big]^{s_j}+v^{-2}F_n^2\big(\vec{h}\big)\bigg].
\end{gather}

\begin{theorem}
\label{th:minimax-abstract}
Let assumptions of Theorem \ref{th:consec-th1-unbounded} be fulfilled and suppose additionally that $K$ satisfies Assumption \ref{ass2:on-kernel-deconvolution-part1}.
Then for any $n\geq 3$, $p>1, \mathbf{q}>1, R>1, D>0,  0<\underline{\Bv}\leq\overline{\Bv}<\infty, \mathbf{u}\in(p/2,\infty], \mathbf{u}\geq \mathbf{q}$, $\vec{s}\in (1,\infty)^d$, $\vec{q}\in [p,\infty)^d$  and any $\bF\subset \bB_{\mathbf{q},d}(D)\cap\bF_{g,\mathbf{u}}(R,D)$
\begin{equation*}
\sup_{f\in\bF}\cR^{(p)}_n[\widehat{f}_{\vec{\mathbf{h}}(\cdot)}, f]\leq C^{(2)}\bigg[l_\bH(\underline{\Bv})+\int_{\underline{\Bv}}^{\overline{\Bv}}v^{p-1}\big[\blL_{\vec{s}}(v,\bF,\mathbf{u})
\wedge\blL_{\vec{s}}(v,\bF)\big]\rd v+\overline{\Bv}^p
\blL_{\vec{q}}(\overline{\Bv},\bF,\mathbf{u})\bigg]^{\frac{1}{p}}+\mathbf{C}_pn^{-\frac{1}{2}}.
\end{equation*}
If additionally $\mathbf{q}\in (p,\infty)$ one has also
\begin{equation*}
\sup_{f\in\bF}\cR^{(p)}_n[\widehat{f}_{\vec{\mathbf{h}}(\cdot)}, f]\leq C^{(2)}\bigg[l_\bH(\underline{\Bv})+\int_{\underline{\Bv}}^{\overline{\Bv}}v^{p-1}\big[\blL_{\vec{s}}(v,\bF,\mathbf{u})
\wedge\blL_{\vec{s}}(v,\bF)\big]\rd v+\overline{\Bv}^{p-\mathbf{q}}
\bigg]^{\frac{1}{p}}+\mathbf{C}_pn^{-\frac{1}{2}}.
\end{equation*}
Moreover, if $\mathbf{q}=\infty$ one has
\begin{equation*}
\sup_{f\in\bF}\cR^{(p)}_n[\widehat{f}_{\vec{\mathbf{h}}(\cdot)}, f]\leq C^{(2)}\bigg[l_\bH(\underline{\Bv})+\int_{\underline{\Bv}}^{\overline{\Bv}}v^{p-1}\big[\blL_{\vec{s}}(v,\bF,\mathbf{u})
\wedge\blL_{\vec{s}}(v,\bF)\big]\rd v+
\blL_{\vec{s}}(\overline{\Bv},\bF,\mathbf{u})\bigg]^{\frac{1}{p}}+\mathbf{C}_pn^{-\frac{1}{2}}.
\end{equation*}
 Finally, if $\bH=\cH^d_{\text{isotr}}$ all the assertions above remain true for any $\vec{s}\in [1,\infty)^d$  if one replaces  in (\ref{eq4:def-blL})--(\ref{eq1:def-blL}) $\mathbf{B}_{j,s_j,\bF}(\cdot)$ by $\mathbf{B}^*_{j,s_j,\bF}(\cdot)$.
\end{theorem}

\noindent It is important to emphasize that $C^{(2)}$ depends only on $\vec{s}, \vec{q},  g, \cK,d$, $R,D, \mathbf{u}$ and $\mathbf{q}$. Note also that the assertions of the theorem remain true if we minimize right hand sides of obtained inequalities w.r.t $\vec{s}, \vec{q}$ since their left hand sides  are independent of $\vec{s}$ and $\vec{q}$. In this context it is important to realize that $C^{(2)}=C^{(2)}(\vec{s},\cdots)$ is bounded for any $\vec{s}\in (1,\infty)^d$ but $C^{(2)}(\vec{s},\cdots)=\infty$ if there exists $j=1,\ldots,d$ such that $s_j=1$. Contrary to that  $C^{(2)}(\vec{s},\cdots)<\infty$ for any $\vec{s}\in [1,\infty)^d$ if $\bH=\cH^d_{\text{isotr}}$ and it explains in particular the fourth  assertion of the theorem.

Note also that  $D,R,\mathbf{u}, \mathbf{q}$ are not involved in the construction of our pointwise selection rule. That means that one and the same estimator can be actually applied on any
$$
\bF\subset\bigcup_{R,D,\mathbf{u},\mathbf{q}}\bB_{\mathbf{q},d}(D)\cap\bF_{g,\mathbf{u}}(R,D).
$$
Moreover, the assertion of the theorem has a non-asymptotical nature; we do not suppose that the number of observations $n$ is large.

\paragraph{\textsf{Discussion}} As we see, the application of our results to some functional class is mainly reduced to the computation of the functions $\mathbf{B}^*_{j,s,\bF}(\cdot)$ $j=1,\ldots,d,$ for some properly chosen $s$. Note however that this task is not necessary for many functional classes used in nonparametric statistics, at least for the classes defined by the help of kernel approximation. Indeed, a typical description of $\bF$ can be summarized as follows.
Let $\lambda_j:\bR_+\to\bR_+$, be such that $\lambda_j(0)=0, \lambda_j\uparrow$ for any $j=1,\ldots,d$. Then, the functional class, say $\bF_K\big[\vec{\lambda}(\cdot),\vec{r}\big]$ can be defined as a collection of functions satisfying
\begin{equation}
\label{eq:class-def-general}
\big\|b_{\mathbf{h},f,j}\big\|_{r_j}\leq \lambda_j(\mathbf{h}),\quad \forall \mathbf{h}\in\cH,
\end{equation}
for some $\vec{r}\in [1,\infty]$.  It yields obviously
$$
\mathbf{B}_{j,r_j,\bF}(\cdot)\leq \lambda_j(\cdot),\quad j=1,\ldots,d,
$$
and the result of Theorem \ref{th:minimax-abstract} remains valid  if we replace formally $\mathbf{B}_{j,r_j,\bF}(\cdot)$ by $\lambda_j(\cdot)$ in all the expressions appearing in this theorem. In Part II we show that for some particular kernel $K^*$, the anisotropic Nikol'skii class $\bN_{\vec{r},d}\big(\vec{\beta},\vec{L}\big)$ is included into the class defined by (\ref{eq:class-def-general}) with $\lambda_j(\mathbf{h})=L_j\mathbf{h}^{\beta_j}$,  whatever the values of $\vec{\beta},\vec{L}$ and $\vec{r}$.

Denote $\vartheta=(\vec{\lambda}(\cdot), \vec{r})$ and remark  that in many cases $\bF_K[\vartheta]\subset\bB_{\mathbf{q},d}(D)$ for any $\vartheta\in\Theta$ for some class parameter $\Theta$ and $\mathbf{q}\geq p, D>0$.
Then, replacing $\mathbf{B}_{j,r_j,\bF}(\cdot)$ by $\lambda_j(\cdot)$ in (\ref{eq4:def-blL}) and (\ref{eq1:def-blL}) and choosing $\vec{q}=(\mathbf{q},\ldots,\mathbf{q})$ we come to the quantities
$
\blL\big(v,\mathbf{u},\vartheta\big)$ and  $\blL_{\mathbf{q}}\big(v,\vartheta\big),
$
 completely determined by the functions $\lambda_j(\cdot), j=1,\ldots,d$, the vector $\vec{r}$ and the number $\mathbf{q}$.
Therefore, putting
$$
\psi_n\big(\vartheta\big)=\inf_{0<\underline{\Bv}\leq\overline{\Bv}<\infty}\bigg(l_\bH(\underline{\Bv})+\int_{\underline{\Bv}}^{\overline{\Bv}}
v^{p-1}\big[\blL(v,\mathbf{u},\theta)\wedge\blL(v,\theta)\big]\rd v+\overline{\Bv}^p
\blL_{\mathbf{q}}\big(\overline{\Bv},\mathbf{u},\vartheta\big)\wedge\overline{\Bv}^{p-\mathbf{q}} \bigg)^{\frac{1}{p}}+n^{-\frac{1}{2}}
$$
we deduce from the first  and the second assertions of Theorem \ref{th:minimax-abstract} for any $\vec{\lambda}(\cdot)$ and $\vec{r}$ and $n\geq 3$
\begin{equation}
\label{eq:discussion-after-th:adaptive-absract-kernel}
\sup_{f\in\bF_K[\vartheta]}\cR^{(p)}_n[\widehat{f}_{\vec{\mathbf{h}}(\cdot)}, f]\leq C^{(3)}\psi_n(\vartheta).
\end{equation}
Since the estimator $\widehat{f}_{\vec{\mathbf{h}}(\cdot)}$ is completely data-driven and, therefore, is independent of $\vec{\lambda}(\cdot)$ and $\vec{r}$,
the bound (\ref{eq:discussion-after-th:adaptive-absract-kernel}) holds for the scale of functional classes $\big\{\bF_K[\vartheta]\big\}_{\vartheta}$.

If $\phi_n\big(\bF_K[\vartheta]\big)$ is the minimax risk defined in (\ref{eq:minmax-risk}) and
\begin{equation}
\label{eq1:discussion-after-th:adaptive-absract-kernel}
\limsup_{n\to\infty}\psi_n(\vartheta)\phi^{-1}_n\big(\bF_{K}[\vartheta]\big)<\infty,\quad \forall
\vartheta\in\Theta,
\end{equation}
 we can assert that our estimator is optimally adaptive over the considered scale $\big\{\bF_K[\vartheta],\;\vartheta\in\Theta \big\}$.

To illustrate the powerfulness of our approach, let us consider a particular scale of functional classes defined by (\ref{eq:class-def-general}).

\paragraph{\textsf{Classes of H\"olderian type}}  Let $\vec{\beta}\in (0,\infty)^{d}$ and $\vec{L}\in (0,\infty)^d$ be given vectors.
\vskip-0.6cm
\begin{definition}
\label{def:class-bfF}
We say that a function $f$ belongs to the class $\bF_{K}\big(\vec{\beta},\vec{L}\big)$, where $K$ satisfies Assumption \ref{ass2:on-kernel-deconvolution-part1}, if $f\in\bB_{\infty,d}\big(max_{j=1,\ldots,d}L_j\big)$ and for any  $j=1,\ldots,d$
\vspace{-0.1cm}
$$
\big\|b_{\mathbf{h},f,j}\big\|_{\infty}\leq L_j\mathbf{v}^{\beta_j},\quad \forall \mathbf{h}\in\cH.
$$
\end{definition}
\vspace{-0.1cm}
We remark that this class is a particular case of the one defined in (\ref{eq:class-def-general}), since it corresponds to $\lambda_j(\mathbf{h})=L_j\mathbf{h}^j$ and $r_j=\infty$ for any $j=\ldots,d$. Moreover let us introduce the following notations
\vspace{-0.3cm}
\begin{eqnarray*}
 \label{eq:beta&omega}\varphi_n=\delta_n^{\frac{1}{2+1/\beta(\alpha)}},\quad \delta_n=L(\alpha)n^{-1}\ln{(n)}, \quad \frac{1}{\beta(\alpha)}=\sum_{j=1}^d\frac{2\blg_j(\alpha)+1}{\beta_j},
\quad L(\alpha)=\prod_{j=1}^dL_j^{\frac{2\blg_j(\alpha)+1}{\beta_j}}.
\end{eqnarray*}

\vspace{-0.3cm}

Then the following result is a direct consequence of Theorem \ref{th:minimax-abstract}. Its simple and short proof is postponed to Section \ref{sec:subcec:proof-of-assertion}.

\begin{assertion}
\label{assert:example-holder} Let the assumptions of Theorem \ref{th:minimax-abstract} be fulfilled. Then
for any $n\ge 3$, $p>1$, $\vec{\beta}\in (0,\infty)^d$, $0<L_0\leq L_\infty<\infty$ and $\vec{L}\in [L_0,L_\infty]^d$ there exists $C>0$ independent of $\vec{L}$ such that
\vskip-0.2cm
\begin{equation*}
\limsup_{n\to\infty}\psi^{-1}_n\big(\vec{\beta}, \vec{L}\big)\sup_{f\in\bF_{K}\left(\vec{\beta},\vec{L}\right)}\cR^{(p)}_n[\widehat{f}_{\vec{\mathbf{h}}(\cdot)}, f]\leq C,
\end{equation*}
\vskip-0.2cm
where we have denoted
\vskip-0.5cm
\begin{eqnarray}
 \label{eq:rate-holder}
\psi_n\big(\vec{\beta},\vec{L}\big)=
\left\{
\begin{array}{ccc}
\ln(n)^{\frac{t(\bH)}{p}}\delta_n^{\frac{(1-1/p)\beta(\alpha)}{\beta(\alpha)+1}},\quad& 2+1/\beta(\alpha)>p;
\\*[1mm]
\ln(n)^{\frac{1\vee t(\bH)}{p}}\delta_n^{\frac{(1-1/p)\beta(\alpha)}{\beta(\alpha)+1}},\quad& 2+1/\beta(\alpha)=p;
\\*[2mm]
\delta_n^{\frac{\beta(\alpha)}{2\beta(\alpha)+1}}, \quad& 2+1/\beta(\alpha)<p.
\end{array}
\right.
\end{eqnarray}
\vskip-0.4cm
\end{assertion}
\vskip-0.6cm
It is interesting to note that the obtained bound, being a very particular case of our consideration in Part II, is completely new if $\alpha\neq0$.
As we already mentioned, for some particular choice of the kernel $K^*$, the anisotropic Nikol'skii class $\bN_{\vec{r},d}\big(\vec{\beta},\vec{L}\big)$ is included in the class $\bF_{K^*}\big[\vec{\lambda}(\cdot),\vec{r}\big]$ with $\lambda_j(\mathbf{v})=L_j\mathbf{v}^{\beta_j}$, whatever the values of $\vec{\beta},\vec{L}$ and $\vec{r}$. Therefore, the aforementioned result
holds on an arbitrary H\"older class $\bN_{\vec{\infty},d}\big(\vec{\beta},\vec{L}\big)$. Comparing the result of Assertion \ref{assert:example-holder} with the lower bound for the minimax risk obtained in \cite{LW15}, we can state that it differs only by some logarithmic factor.
Using the modern statistical language, we say that the estimator $\widehat{f}_{\vec{\mathbf{h}}(\cdot)}$ is {\it nearly optimally-adaptive} over the scale of H\"older classes.

\section{Proofs}

\subsection{Proof of Theorem \ref{th:L_p-norm-oracle-deconv}}
\label{sec:subsec-proof-l_p-norm-ineq-deconv}

The main ingredients of the proof of the theorem are given in Proposition \ref{prop:aux-result2}. Their proofs are postponed to Section \ref{sec:subsec-proof-prop:aux-result2}.
Introduce  for any  $\vec{h}\in\cH^d$
\vskip-0.5cm
\begin{eqnarray*}
\xi_{n}\big(x,\vec{h}\big)&=&\frac{1}{n}\sum_{i=1}^n\big[M\big(Z_i-x,\vec{h}\big)-\bE_f M\big(Z_i-x,\vec{h}\big)\big],\quad  x\in\bR^d.
\end{eqnarray*}

\begin{proposition}
\label{prop:aux-result2}
Let Assumptions \ref{ass1:ass-on-noise-upper-bound} and \ref{ass:on-kernel-deconvolution} be fulfilled. Then for any $n\geq 3$ and any $p> 1$
\begin{eqnarray*}
&&(\mathbf{i})\quad\int_{\bR^d}\bE_f\Big\{\sup_{\vec{h}\in\cH^d}\big[\big|\xi_{n}\big(x,\vec{h}\big)\big|-U_n\big(x,\vec{h}\big)\big]_+^p\Big\}
\nu_d(\rd x)\leq C_p n^{-\frac{p}{2}};
\end{eqnarray*}
\begin{eqnarray*}
&&(\mathbf{ii})\quad \int_{\bR^d}\bE_f\Big\{\sup_{\vec{h}\in\cH^d}\big[\widehat{U}_n\big(x,\vec{h}\big)-3U_n\big(x,\vec{h}\big)\big]_+^p\Big\}\nu_d(\rd x)\leq
C^{\prime}_p n^{-\frac{p}{2}};
\\*[2mm]
&&(\mathbf{iii})\quad \int_{\bR^d}\bE_f\Big\{\sup_{\vec{h}\in\cH^d}\big[U_n\big(x,\vec{h}\big)-4\widehat{U}_n\big(x,\vec{h}\big)\big]_+^p\Big\}\nu_d(\rd x)\leq C^{\prime}_p n^{-\frac{p}{2}}.
\end{eqnarray*}

\end{proposition}
The explicit expression of constant $C_p$ and $C^{\prime}_p$ can be found in the proof.

\subsubsection{Proof of the theorem}

We start by proving the so-called pointwise oracle inequality.

\vskip0.1cm

\noindent {\textsf{Pointwise oracle inequality.\;}
Let $\vec{h}\in\bH$  and $x\in\bR^d$ be fixed.
We have  in view of the triangle inequality
\begin{equation}
\label{eq1:proof-th:L_p-norm-oracle-inequality-deconv}
\left|\widehat{f}_{\vec{\mathbf{h}}(x)}(x)-f(x)\right|\leq\left|\widehat{f}_{\vec{\mathbf{h}}(x)\vee\vec{h}}(x)-\widehat{f}_{\vec{\mathbf{h}}(x)}(x)\right|+
\left|\widehat{f}_{\vec{\mathbf{h}}(x)\vee\vec{h}}(x)-\widehat{f}_{\vec{h}}(x)\right|+\left|\widehat{f}_{\vec{h}}(x)-f(x)\right|.
\end{equation}

$1^0.\;$ First, note that obviously  $\widehat{f}_{\vec{\mathbf{h}}(x)\vee\vec{h}}(x)= \widehat{f}_{\vec{h}\vee\vec{\mathbf{h}}(x)}(x)$ and, therefore,
\begin{eqnarray*}
\left|\widehat{f}_{\vec{\mathbf{h}}(x)\vee\vec{h}}(x)-\widehat{f}_{\vec{\mathbf{h}}(x)}(x)\right|=
\left|\widehat{f}_{\vec{h}\vee\vec{\mathbf{h}}(x)}(x)-\widehat{f}_{\vec{\mathbf{h}}(x)}(x)\right|\leq \widehat{\cR}_{\vec{h}}(x)+ 4\widehat{U}_n\big(x,\vec{\mathbf{h}}(x)\vee\vec{h}\big)+4\widehat{U}_n\big(x,\vec{\mathbf{h}}(x)\big).
\end{eqnarray*}
Moreover by definition, $\widehat{U}_n\big(x,\vec{\eta}\big)\leq \widehat{U}^*_n\big(x,\vec{\eta}\big)$ for any $\vec{\eta}\in\cH^d$.

Next, for any $\vec{h},\vec{\eta}\in\cH^d$ we have obviously
$
\widehat{U}_n\big(x,\vec{h}\vee\vec{\eta}\big)\leq \widehat{U}^*_n\big(x,\vec{h}\big)\wedge \widehat{U}^*_n\big(x,\vec{\eta}\big).
$
Thus, we obtain
\begin{eqnarray}
\label{eq2:proof-th:L_p-norm-oracle-inequality-deconv}
\left|\widehat{f}_{\vec{\mathbf{h}}(x)\vee\vec{h}}(x)-\widehat{f}_{\vec{\mathbf{h}}(x)}(x)\right|\leq \widehat{\cR}_{\vec{h}}(x)+ 8\widehat{U}^*_n\big(x,\vec{\mathbf{h}}(x)\big).
\end{eqnarray}
Similarly we have
\begin{eqnarray}
\label{eq3:proof-th:L_p-norm-oracle-inequality-deconv}
&&\left|\widehat{f}_{\vec{\mathbf{h}}(x)\vee\vec{h}}(x)-\widehat{f}_{\vec{h}}(x)\right|\leq \widehat{\cR}_{\vec{\mathbf{h}}(x)}(x)+
8\widehat{U}^*_n\big(x,\vec{h}\big).
\end{eqnarray}
The definition of $\vec{\mathbf{h}}(x)$ implies that for any $\vec{h}\in\bH$
$$
\widehat{\cR}_{\vec{\mathbf{h}}(x)}(x)+ 8\widehat{U}^*_n\big(x,\vec{\mathbf{h}}(x)\big)+\widehat{\cR}_{\vec{h}}(x)+
8\widehat{U}^*_n\big(x,\vec{h}\big)\leq 2\widehat{\cR}_{\vec{h}}(x)+ 16\widehat{U}^*_n\big(x,\vec{h}\big)
$$
and we get from (\ref{eq1:proof-th:L_p-norm-oracle-inequality-deconv}), (\ref{eq2:proof-th:L_p-norm-oracle-inequality-deconv}) and (\ref{eq3:proof-th:L_p-norm-oracle-inequality-deconv}) for any $\vec{h}\in\bH$
\begin{equation}
\label{eq4:proof-th:L_p-norm-oracle-inequality-deconv}
\left|\widehat{f}_{\vec{\mathbf{h}}(x)}(x)-f(x)\right|\leq 2\widehat{\cR}_{\vec{h}}(x)+ 16\widehat{U}^*_n\big(x,\vec{h}\big)+\left|\widehat{f}_{\vec{h}}(x)-f(x)\right| .
\end{equation}

$2^0.\;$ We obviously have for any $\vec{h},\vec{\eta}\in\cH^d$
$$
\left|\widehat{f}_{\vec{h}\vee\vec{\eta}}(x)-\widehat{f}_{\vec{\eta}}(x)\right|\leq \big|\bE_f M\big(Z_1-x,\vec{h}\vee\vec{\eta}\big)-\bE_f M\big(Z_1-x,\vec{\eta}\big)\big| +\big|\xi_n\big(x,\vec{h}\vee\vec{\eta}\big)\big|+\big|\xi_n\big(x,\vec{\eta}\big)\big|.
$$
Note that for any $\mathrm{h}\in\cH^d$
\begin{eqnarray*}
\bE_f M\big(Z_1-x,\vec{\mathrm{h}}\big)&:=&\int_{\bR^d}M\big(t-x,\vec{\mathrm{h}}\big) \mathfrak{p}(t)\nu_d(\rd t)
\\
\nonumber
&=&(1-\alpha)\int_{\bR^d}M\big(t-x,\vec{\mathrm{h}}\big)f(t)\nu_d(\rd t)+\alpha\int_{\bR^d}M\big(t-x,\vec{\mathrm{h}}\big)\big[f\star g\big](t)\nu_d(\rd t),
\end{eqnarray*}
in view of the structural assumption (\ref{eq:convolution structure model}) imposed on the density $\mathfrak{p}$. Note that
\begin{eqnarray*}
&&(1-\alpha)\int_{\bR^d}M\big(t-x,\vec{\mathrm{h}}\big)f(t)\nu_d(\rd t)+\alpha\int_{\bR^d}M\big(t-x,\vec{\mathrm{h}}\big)\big[f\star g\big](t)\nu_d(\rd t)
\\
&&=\int_{\bR^d}f(z)
\bigg[(1-\alpha)M\big(z-x,\vec{\mathrm{h}}\big)+\alpha\int_{\bR^d}M\big(u,\vec{\mathrm{h}}\big)g(u-[z-x])\nu_d(\rd u)\bigg]\nu_d(\rd z)
\end{eqnarray*}
and, therefore, in view of the definition of $M\big(\cdot,\vec{h}\big)$, c.f. (\ref{eq:def-kernel-M}), we obtain for any $\mathrm{h}\in\cH^d$
\begin{eqnarray}
\label{eq5:proof-th:L_p-norm-oracle-inequality-deconv}
\bE_f M\big(Z_1-x,\vec{\mathrm{h}}\big)=\int_{\bR^d}K_{\vec{\mathrm{h}}}(z-x)f(z)\nu_d(\rd z)=:S_{\vec{\mathrm{h}}}(x,f).&&
\end{eqnarray}
We deduce from (\ref{eq5:proof-th:L_p-norm-oracle-inequality-deconv}) that
$$
\big|\bE_f M\big(Z_1-x,\vec{h}\vee\vec{\eta}\big)-\bE_f M\big(Z_1-x,\vec{\eta}\big)\big|=\big|S_{\vec{h}\vee \vec{\eta}}(x,f)-S_{ \vec{\eta}}(x,f)\big|
$$
and, therefore,
for any $\vec{h},\vec{\eta}\in\cH^d$
\begin{eqnarray}
\label{eq6:proof-th:L_p-norm-oracle-inequality-deconv}
\left|\widehat{f}_{\vec{h}\vee\vec{\eta}}(x)-\widehat{f}_{\vec{\eta}}(x)\right|\leq \big|S_{\vec{h}\vee \vec{\eta}}(x,f)-S_{ \vec{\eta}}(x,f)\big|+\big|\xi_n\big(x,\vec{h}\vee\vec{\eta}\big)\big|+\big|\xi_n\big(x,\vec{\eta}\big)\big|.
\end{eqnarray}

$3^0.\;$ Set for any $\vec{h}\in\cH^d$ and any $x\in\bR^d$
\begin{gather*}
\upsilon(x)=\sup_{\vec{\eta}\in\cH^d}\big[\big|\xi_n\big(x,\vec{\eta}\big)\big|-U_n\big(x,\vec{\eta}\big)\big]_+
\\
\varpi_1(x)=\sup_{\vec{h}\in\cH^d}\big[U_n\big(x,\vec{h}\big)-4\widehat{U}_n\big(x,\vec{h}\big)\big]_+,\quad \varpi_2(x)=\sup_{\vec{h}\in\cH^d}\big[\widehat{U}_n\big(x,\vec{h}\big)-3U_n\big(x,\vec{h}\big)\big]_+
\end{gather*}
We obtain in view of (\ref{eq6:proof-th:L_p-norm-oracle-inequality-deconv}) that for any $\vec{h}\in\bH$ (since obviously $\vec{h}\vee\vec{\eta}\in\cH^d$ for any $\vec{h},\vec{\eta}\in\cH^d$)
\begin{eqnarray}
\label{eq7:proof-th:L_p-norm-oracle-inequality-deconv}
\widehat{\cR}_{\vec{h}}(x)\leq B^*_{\vec{h}}(x,f)+2\upsilon(x)+2\varpi_1(x).&&
\end{eqnarray}
Note also that in view of the obvious inequality $(\sup_{\alpha}F_\alpha-\sup_{\alpha}G_\alpha)_+\leq \sup_{\alpha}(F_\alpha-G_\alpha)_+$
\begin{eqnarray}
\label{eq8:proof-th:L_p-norm-oracle-inequality-deconv}
\big[\widehat{U}^*_n\big(x,\vec{h}\big)-3U^*_n\big(x,\vec{h}\big)\big]_+\leq \sup_{\vec{\eta}\in\cH^d}\big[\widehat{U}_n\big(x,\vec{\eta}\big)-3U_n\big(x,\vec{\eta}\big)\big]_+=:\varpi_2(x)
\end{eqnarray}
We get from (\ref{eq4:proof-th:L_p-norm-oracle-inequality-deconv}), (\ref{eq7:proof-th:L_p-norm-oracle-inequality-deconv}) and (\ref{eq8:proof-th:L_p-norm-oracle-inequality-deconv})
\begin{equation*}
\left|\widehat{f}_{\vec{\mathbf{h}}(x)}(x)-f(x)\right|\leq 2B^*_{\vec{h}}(x,f)+4\upsilon(x)+4\varpi_1(x)+48U^*_n\big(x,\vec{h}\big)+ 16\varpi_2(x)+\left|\widehat{f}_{\vec{h}}(x)-f(x)\right|.
\end{equation*}
It remains to note that
$$
\left|\widehat{f}_{\vec{h}}(x)-f(x)\right|\leq B_{\vec{h}}(x,f)+\big|\xi_n\big(x,\vec{h}\big)\big|\leq B_{\vec{h}}(x,f)+U_n\big(x,\vec{h}\big)+\upsilon(x),
$$
and we obtain for any $\vec{h}\in\bH$ and  $x\in\bR^d$
\begin{equation*}
\left|\widehat{f}_{\vec{\mathbf{h}}(x)}(x)-f(x)\right|\leq 2B^*_{\vec{h}}(x,f)+B_{\vec{h}}(x,f)+5\upsilon(x)+4\varpi_1(x)+49U^*_n\big(x,\vec{h}\big)+ 16\varpi_2(x).
\end{equation*}
Noting that the left hand side of the latter inequality is independent of $\vec{h}$ we obtain for any $x\in\bR^d$
\begin{equation}
\label{eq9:proof-th:L_p-norm-oracle-inequality-deconv}
\left|\widehat{f}_{\vec{\mathbf{h}}(x)}(x)-f(x)\right|\leq
\inf_{\vec{h}\in\bH}\Big\{2B^*_{\vec{h}}(x,f)+B_{\vec{h}}(x,f)+49U^*_n\big(x,\vec{h}\big)\Big\}
+5\upsilon(x)+4\varpi_1(x)+ 16\varpi_2(x).
\end{equation}
This is the pointwise oracle inequality.

\vskip0.1cm

\noindent {\textsf{Application of Proposition \ref{prop:aux-result2}.}\; Set for any $x\in\bR^d$
$$
R_n(x)=\inf_{\vec{h}\in\bH}\Big\{2B^*_{\vec{h}}(x,f)+B_{\vec{h}}(x,f)+49U^*_n\big(x,\vec{h}\big)\Big\}
$$
Applying  Proposition \ref{prop:aux-result2}  we obtain in view of (\ref{eq9:proof-th:L_p-norm-oracle-inequality-deconv}) and the triangle inequality
 \begin{eqnarray*}
\label{eq10:proof-th:L_p-norm-oracle-inequality-deconv}
\cR^{(p)}_n[\widehat{f}_{\vec{\mathbf{h}}(\cdot)}, f]&\leq & \big\|R_n\big\|_p+5\bigg[\int_{\bR^d}\bE_f\big\{\upsilon(x)\big\}^p\bigg]^{\frac{1}{p}}
+4\bigg[\int_{\bR^d}\bE_f\big\{\varpi_1(x)\big\}^p\bigg]^{\frac{1}{p}}+16\bigg[\int_{\bR^d}\bE_f\big\{\varpi_2(x)\big\}^p\bigg]^{\frac{1}{p}}
\\
&\leq &\big\|R_n\big\|_p+\mathbf{C}_pn^{-\frac{1}{2}},
\end{eqnarray*}
where $\mathbf{C}_p=5(C_p)^{\frac{1}{p}}+20(C^\prime_p)^{\frac{1}{p}}$. The theorem is proved.
\epr

\subsubsection{Proof of Proposition \ref{prop:aux-result2}}
\label{sec:subsec-proof-prop:aux-result2}

Since the proof of the proposition is quite long and technical, we divide it into several steps.

\smallskip

\paragraph{Preliminaries}

$1^0.\;$ We start the proof with the following simple remark. Let $\check{M}\big(t,\vec{h}\big), t\in\bR^d,$ denote the Fourier transform of $M\big(\cdot,\vec{h}\big)$.
Then, we obtain in view of the definition of $M\big(\cdot,\vec{h}\big)$
\begin{eqnarray}
\label{eq:fourier-trans-of-kernel-M}
\check{M}\big(t,\vec{h}\big)=\check{K}\big(t\vec{h}\big)\big[(1-\alpha)+\alpha\check{g}(-t)\big]^{-1},\;\; t\in\bR^d.&&
\end{eqnarray}
Note that   Assumptions \ref{ass1:ass-on-noise-upper-bound} and \ref{ass:on-kernel-deconvolution} guarantee that $\check{M}\big(\cdot,\vec{h}\big)\in\bL_1\big(\bR^d\big)\cap\bL_2\big(\bR^d\big)$ for any $\vec{h}\in\cH^d$  and, therefore,
$$
\big\|M\big(\cdot,\vec{h}\big)\big\|_\infty\leq (2\pi)^{-d}\big\|\check{M}\big(\cdot,\vec{h}\big)\big\|_1,\qquad
\big\|M\big(\cdot,\vec{h}\big)\big\|_2= (2\pi)^{-d}\big\|\check{M}\big(\cdot,\vec{h}\big)\big\|_2.
$$
Thus, putting
$$
\cM_\infty(\vec{h}\big)= M_\infty\prod_{j=1}^dh_j^{-1}(h_j\wedge 1)^{-\blg_j(\alpha)},
$$
 we obtain in view of Assumptions \ref{ass1:ass-on-noise-upper-bound} and \ref{ass:on-kernel-deconvolution} for any $\vec{h}\in\cH^d$
\begin{eqnarray}
\label{eq:bound-for-supnorm-and-L_2-norm}
\big\|M\big(\cdot,\vec{h}\big)\big\|_\infty\leq \cM_\infty(\vec{h}\big),\qquad
\big\|M\big(\cdot,\vec{h}\big)\big\|_2\leq M_2\prod_{j=1}^dh_j^{-\frac{1}{2}}(h_j\wedge 1)^{-\blg_j(\alpha)},
&&
\end{eqnarray}
where
$
M_2=\big[(2\pi)^{-d}\big\{\e^{-1}\big\|\check{K}\big\|_2\mathrm{1}_{\alpha\neq 1}+\Upsilon_0^{-1}\mathbf{k}_2\mathrm{1}_{\alpha=1}\big\}\big]\vee 1.
$
Additionally we deduce from (\ref{eq:bound-for-supnorm-and-L_2-norm})
\begin{eqnarray}
\label{eq:bound-for-L_4-norm}
&&\big\|M\big(\cdot,\vec{h}\big)\big\|^4_4\leq \big\|M\big(\cdot,\vec{h}\big)\big\|^2_\infty\big\|M\big(\cdot,\vec{h}\big)\big\|^2_2
\leq M^2_2M^2_\infty \prod_{j=1}^dh_j^{-3}(h_j\wedge 1)^{-4\blg_j(\alpha)},\quad \forall\vec{h}\in\cH^d.
\end{eqnarray}

Let $\cL\big(\cdot,\vec{h}\big)$ be either $M\big(\cdot,\vec{h}\big)$ or $M^2\big(\cdot,\vec{h}\big)$ and
let $\cL_\infty\big(\vec{h}\big)$ denote either $\cM_\infty\big(\vec{h}\big)$ or $\cM^2_\infty\big(\vec{h}\big)$.

\noindent We have in view of (\ref{eq:bound-for-supnorm-and-L_2-norm})
\begin{eqnarray}
\label{eq5900:proof-th:aux-result1}
\cL^{-1}_\infty\big(\vec{h}\big)\vee\cL_\infty\big(\vec{h}\big)&\leq& M_\infty^2e^{2\sum_{j=1}^d(1+\blg_j(\alpha))|\ln(h_j)|},\;\;\;
\forall \vec{h}\in\cH^d.
\end{eqnarray}
Additionally, we get from (\ref{eq:bound-for-supnorm-and-L_2-norm}) and (\ref{eq:bound-for-L_4-norm})
\begin{eqnarray}
\label{eq10000:proof-th:aux-result1}
\big\|\cL\big(\cdot,\vec{h}\big)\big\|_2^2\leq M_2^2M_\infty^2e^{4\sum_{j=1}^d(1+\blg_j(\alpha))|\ln(h_j)|} ,\quad
\forall \vec{h}\in\cH^d.
\end{eqnarray}
Set
$
\sigma^{\cL}\big(x,\vec{h}\big)=\sqrt{\int_{\bR^d}\cL^2\big(t-x,\vec{h}\big) \mathfrak{p}(t)\nu_d(\rd t)}
$
and  note that in view of (\ref{eq10000:proof-th:aux-result1})  for any $\vec{h}\in\cH^d$
\begin{equation}
\label{eq2:proof-th:aux-result1}
\int_{\bR^d}\big[\sigma^{\cL}\big(x,\vec{h}\big)\big]^2\nu_d(\rd x)=\big\|\cL\big(\cdot,\vec{h}\big)\big\|_2^2\leq M_2^2M_\infty^2e^{4\sum_{j=1}^d(1+\blg_j(\alpha))|\ln(h_j)|} .
\end{equation}
Next, we have in view of (\ref{eq5900:proof-th:aux-result1})
\begin{equation}
\label{eq0002:proof-th:aux-result1}
\big\|\sigma^{\cL}\big(\cdot,\vec{h}\big)\big\|_\infty\leq \cL_\infty\big(\vec{h}\big)\leq M_\infty^2e^{2\sum_{j=1}^d(1+\blg_j(\alpha))|\ln(h_j)|} .
\end{equation}

$2^0.\;$
Define for any $x\in\bR^d$ and $\vec{h}\in\cH^d$
\vskip-0.7cm
\begin{gather*}
\xi^{\cL}\big(x,\vec{h}\big)=n^{-1}\sum_{i=1}^n\big[\cL\big(Z_i-x,\vec{h}\big)-\bE \cL\big(Z_i-x,\vec{h}\big)\big];
\\*[-3mm]
z_n\big(x,\vec{h}\big)=3\ln(n)+(8p+22)\sum_{j=1}^d\big[1+\blg_j(\alpha)\big]\big|\ln(\mathrm{h}_j)\big|+
2\big|\ln{\big(\big\{\sigma^{\cL}\big(x,\vec{h}\big)\big\}\vee \big\{n^{-3/2}\cL_\infty\big(\vec{h}\big)\big\}\big)}\big|.
\\[-1mm]
V^{\cL}\big(x,\vec{h}\big)=\sigma^{\cL}\big(x,\vec{h}\big)\sqrt{\frac{2z_n\big(x,\vec{h}\big)}{n}}+\frac{4z_n\big(x,\vec{h}\big)
\cL_\infty\big(\vec{h}\big)}{3n};
\\[-0mm]
\hskip-0.6cm U^{\cL}\big(x,\vec{h}\big)=\sigma^{\cL}\big(x,\vec{h}\big)\sqrt{\frac{2\lambda_n\big(\vec{h}\big)}{n}}+
\frac{4\lambda_n\big(\vec{h}\big)\cL_\infty\big(\vec{h}\big)}{3n},
\end{gather*}
where remind
$
\lambda_n\big(\vec{h}\big)=4\ln(M_\infty)+6\ln{(n)}+(8p+26)\sum_{j=1}^d\big[1+\blg_j(\alpha)\big]\big|\ln(h_j)\big|.
$

\smallskip

Noting that $\sup_{z\in[a,b]}|\ln z|\leq |\ln a|\vee|\ln b|$ for any $0<a<b<\infty$ we deduce from (\ref{eq0002:proof-th:aux-result1})
$
z_n\big(x,\vec{h}\big)\leq \lambda_n\big(\vec{h}\big)
$
for any $x\in\bR^d$  and, therefore, for any $\vec{h}\in\cH^d$
\begin{eqnarray}
\label{eq1:proof-aux_proposition1}
V^{\cL}\big(x,\vec{h}\big)\leq U^{\cL}\big(x,\vec{h}\big).
\end{eqnarray}

\paragraph{\textsf{First step}}
Let $x\in\bR^d$ and $\vec{h}\in\cH^d$ be fixed and put $b=8p+22$.

We obtain for any $z\geq 1$ and $q\geq 1$ by the integration of the Bernstein inequality
\begin{equation*}
\label{eq2:discussion-price-to-pay}
\bE_f\bigg\{\big|\zeta^{\cL}\big(x,\vec{h}\big)\big|- \frac{\sqrt{2z}\sigma^\cL\big(x,\vec{h}\big)}{\sqrt{n}}-\frac{4z\cL_\infty\big(\vec{h}\big)}{3n}\bigg\}_+^{q}\leq 2\Gamma(q+1)\bigg[\frac{\sqrt{2}\sigma^\cL\big(x,\vec{h}\big)}{\sqrt{n}}+\frac{4\cL_\infty\big(\vec{h}\big)}{3n}\bigg]^{q}\exp{\left\{-z\right\}},
\end{equation*}
where $\Gamma$ is the Gamma-function.

\smallskip

$1^0.\;$ Choose
$
z=z_n\big(x,\vec{h}\big).
$
Noting that for any  $n\in \bN^*$ and $x\in\bR^d$
$$
\frac{\sqrt{2}\sigma^\cL\big(x,\vec{h}\big)}{\sqrt{n}}+\frac{4\cL_\infty\big(\vec{h}\big)}{3n}\leq 3\cL_\infty\big(\vec{h}\big)n^{-\frac{1}{2}}
$$
and taking into account that $\exp{\{-|\ln(y)|\}}\leq y$ for any $y>0$, we get
\begin{eqnarray}
\label{eq1:proof-th:aux-result1}
&&\bE_f\big\{\big|\zeta^{\cL}\big(x,\vec{h}\big)\big|- V^\cL\big(x,\vec{h}\big)\big\}_+^{q}\nonumber\\*[2mm]
&&\leq 2\times 3^q\Gamma(q+1)n^{-\frac{q}{2}-3}\cL^{q}_\infty\big(\vec{h}\big)e^{b\sum_{j=1}^d(1+\blg_j(\alpha))|\ln(h_j)|}\big(\big\{\sigma^{\cL}\big(x,\vec{h}\big)\big\}\vee \big\{n^{-3/2}\cL_\infty\big(\vec{h}\big)\big\}\big)^2
\nonumber\\*[2mm]
&&\leq C^{(1)}_qn^{-\frac{q}{2}-3}e^{(2q-b)\sum_{j=1}^d(1+\blg_j(\alpha))|\ln(h_j)|}\big(\big\{\sigma^{\cL}\big(x,\vec{h}\big)\big\}\vee \big\{n^{-3/2}\cL_\infty\big(\vec{h}\big)\big\}\big)^2.
\end{eqnarray}
Here to get the second inequality we have used (\ref{eq5900:proof-th:aux-result1}) and put $C^{(1)}_q=2M_\infty^{2q}3^q\Gamma(q+1)$.

\smallskip

 Set $\cX\big(\vec{h}\big)=\big\{x\in\bR^d:\; \sigma^{\cL}\big(x,\vec{h}\big)\geq n^{-3/2}\cL_\infty\big(\vec{h}\big)\big\}$,  $\bar{\cX}\big(\vec{h}\big)=\bR^d\setminus\cX\big(\vec{h}\big)$ and later on the integration over the empty set is supposed to be zero.

We have in view of (\ref{eq1:proof-aux_proposition1}), (\ref{eq2:proof-th:aux-result1}) and  (\ref{eq1:proof-th:aux-result1}) applied with $q=p$ that for any
$\vec{h}\in\cH^d$
\begin{equation}
\label{eq3:proof-th:aux-result1}
\int_{\cX\big(\vec{h}\big)}\bE_f\big\{\big|\zeta^{\cL}\big(x,\vec{h}\big)\big|-U^\cL\big(x,\vec{h}\big)\big\}_+^p\nu_d(\rd x)
\leq C^{(2)}_pn^{-\frac{p}{2}}e^{(2p+4-b)\sum_{j=1}^d(1+\blg_j(\alpha))|\ln(h_j)|}.
\end{equation}
where $C^{(2)}_p=C^{(1)}_pM^2_2M_\infty^{2}$.

\smallskip

$2^0.\;$ Introduce the following notations. For any $i=1,\ldots, n$ set
$$
\Psi_i\big(x,\vec{h}\big)=\mathrm{1}\Big\{\big|\cL\big(Z_i-x,\vec{h}\big)-\bE\cL\big(Z_i-x,\vec{h}\big)\big|\geq n^{-1}\cL_\infty\big(\vec{h}\big)\Big\},
$$
and introduce the random event
$
D\big(x,\vec{h}\big)=\Big\{\sum_{i=1}^n\Psi_i\big(x,\vec{h}\big)\geq 2\Big\}
$.
As usual, the complimentary event will be denoted by $\bar{D}\big(x,\vec{h}\big)$. Set finally $\pi\big(x,\vec{h}\big)=\bP_f\big\{\Psi_1\big(x,\vec{h}\big)=1\big\}$.

 We obviously have
$$
\big|\zeta^{\cL}\big(x,\vec{h}\big)\big|\mathrm{1}_{\bar{D}\big(x,\vec{h}\big)}\leq \frac{3\cL_\infty\big(\vec{h}\big)}{n}<U^\cL\big(\vec{h}\big)
$$
and, therefore,
\begin{equation}
\label{eq4:proof-th:aux-result1}
\mathrm{1}_{\bar{D}\big(x,\vec{h}\big)}\big\{\big|\zeta^{\cL}\big(x,\vec{h}\big)\big|-U^\cL\big(x,\vec{h}\big)\big\}^p_+=0.
\end{equation}
Applying Cauchy-Schwartz inequality, we deduce from (\ref{eq4:proof-th:aux-result1}) that
\begin{equation*}
\bE_f\big\{\big|\zeta^{\cL}\big(x,\vec{h}\big)\big|-U^\cL\big(x,\vec{h}\big)\big\}^p_+\leq
\Big[\bE_f\big\{\big|\zeta^{\cL}\big(x,\vec{h}\big)\big|-U^\cL\big(x,\vec{h}\big)\big\}^{2p}_+\;\bP_f\big\{D\big(x,\vec{h}\big)\big\}\Big]^{\frac{1}{2}}.
\end{equation*}
Using (\ref{eq1:proof-th:aux-result1}) with $q=2p$ and (\ref{eq5900:proof-th:aux-result1}) we obtain for any $x\in\bar{\cX}\big(\vec{h}\big)$
\begin{equation}
\label{eq40:proof-th:aux-result1}
\bE_f\big\{\big|\zeta^{\cL}\big(x,\vec{h}\big)\big|-U^\cL\big(x,\vec{h}\big)\big\}^p_+\leq
C^{(3)}_pn^{-\frac{p}{2}-3}e^{(2p+2-b/2)\sum_{j=1}^d(1+\blg_j(\alpha))|\ln(h_j)|}\Big[\bP_f\big\{D\big(x,\vec{h}\big)\big\}\Big]^{\frac{1}{2}},
\end{equation}
where we have put $C^{(3)}_p=\big[C^{(1)}_{2p}\big]^{\frac{1}{2}}M_\infty^2$.

For any $\lambda>0$ we have in view of the exponential  Markov inequality
\begin{eqnarray*}
\bP_f\big\{D\big(x,\vec{h}\big)\big\} &=& \bP_f\Big\{\sum_{i=1}^n\Psi_i\big(x,\vec{h}\big)\geq 2\Big\}
\leq e^{-2\lambda} \big[e^{\lambda}\pi\big(x,\vec{h}\big)  + 1-\pi\big(x,\vec{h}\big)\big]^n
\\
&=&
e^{-2\lambda}\big[(e^\lambda-1)\pi\big(x,\vec{h}\big)+1\big]^n \leq \exp\{-2\lambda + n(e^{\lambda}-1)\pi\big(x,\vec{h}\big)\}.
\end{eqnarray*}
We get applying the Tchebychev inequality
$
\pi\big(x,\vec{h}\big)\leq n^2\cL^{-2}_\infty\big(\vec{h}\big)\big[\sigma^{\cL}\big(x,\vec{h}\big)\big]^2.
$
It yields
\begin{equation*}
\label{eq5:proof-th:aux-result1}
\bP_f\big\{D\big(x,\vec{h}\big)\big\}\leq \exp{\big\{-2\lambda + n^{3}\cL^{-2}_\infty\big(\vec{h}\big)\big[\sigma^{\cL}\big(x,\vec{h}\big)\big]^2(e^{\lambda}-1)\big\}}
,\quad \forall \vec{h}\in\cH^d.
\end{equation*}
Note that the definition of $\bar{\cX}\big(\vec{h}\big)$ implies  $n^{3}\cL^{-2}_\infty\big(\vec{h}\big)\big[\sigma^{\cL}\big(x,\vec{h}\big)\big]^2<1$ for any  $x\in\bar{\cX}\big(\vec{h}\big)$.
Hence, choosing  $\lambda=\ln 2 -2\ln{\big\{n^{3/2}\cL^{-1}_\infty\big(\vec{h}\big)\sigma^{\cL}\big(x,\vec{h}\big)\big\}}$ we have
\begin{equation*}
\label{eq5000:proof-th:aux-result1*}
 \bP_f\big\{D\big(x,\vec{h}\big)\big\} \leq
(e^2/4)n^{6}\cL^{-4}_\infty\big[\sigma^{\cL}\big(x,\vec{h}\big)\big]^4,\quad \forall x\in\bar{\cX}\big(\vec{h}\big).
\end{equation*}
It yields, together with (\ref{eq5900:proof-th:aux-result1}), (\ref{eq2:proof-th:aux-result1}) and   (\ref{eq40:proof-th:aux-result1}) and  for any $\vec{h}\in\cH^d$
\begin{eqnarray}
\label{eq6-new:proof-th:aux-result1}
&&\int_{\bar{\cX}\big(\vec{h}\big)}\bE_f\big\{\big|\zeta^{\cL}\big(x,\vec{h}\big)\big|-U^\cL\big(x,\vec{h}\big)\big\}_+^p\nu_d(\rd x)\leq
C^{(4)}_pn^{-\frac{p}{2}}e^{(2p+10-b/2)\sum_{j=1}^d(1+\blg_j(\alpha))|\ln(h_j)|}
\end{eqnarray}
where $C^{(4)}_p=C^{(3)}_p(e/2)M_\infty^{6}M_2^2$.
Putting $C^{(5)}_p=C^{(2)}_p + C^{(4)}_p$ and noting that $2p+10-b/2<0$ we obtain from
(\ref{eq3:proof-th:aux-result1}) and   (\ref{eq6-new:proof-th:aux-result1}) for any
$\vec{h}\in\cH^d$
\begin{eqnarray}
\label{eq6-new-prime:proof-th:aux-result1}
&&\int_{\bR^d}\bE_f\big\{\big|\zeta^{\cL}\big(x,\vec{h}\big)\big|-U^\cL\big(x,\vec{h}\big)\big\}_+^p\nu_d(\rd x)\leq
C^{(5)}_pn^{-\frac{p}{2}}e^{(2p+10-b/2)\sum_{j=1}^d(1+\blg_j(\alpha))|\ln(h_j)|}.
\end{eqnarray}

\smallskip

$3^0.\;$ Choosing $\cL=M$ and  $\cL_\infty=\cM_\infty$  we get from (\ref{eq6-new-prime:proof-th:aux-result1}) and the definition of $b$
\begin{eqnarray}
\label{eq300-new:proof-th:aux-result1}
&&\int_{\bR^d}\bE_f\big\{\big|\xi_{n}\big(x,\vec{h}\big)\big|-U_n\big(x,\vec{h}\big)\big\}_+^p\nu_d(\rd x)
\leq C^{(5)}_p n^{-\frac{p}{2}}e^{-\sum_{j=1}^d(1+\blg_j(\alpha))|\ln(h_j)|},\;\;\;
\forall
\vec{h}\in\cH^d.
\end{eqnarray}
The first assertion of the proposition follows from (\ref{eq300-new:proof-th:aux-result1}) with
$
C_p=C^{(5)}_p\sum_{k\in\bZ^d}e^{-\sum_{j=1}^d|k_j|}.
$

\paragraph{\textsf{Second step}}

Denoting $\chi\big(x,\vec{h}\big)=\big\{\big|\widehat{\sigma}^2\big(x,\vec{h}\big)-\sigma^2\big(x,\vec{h}\big)\big|-\mU_n\big(x,\vec{h}\big)\big\}_+$, where
\begin{eqnarray*}
\mU_n\big(x,\vec{h}\big)&=&\sigma^{M^2}\big(x,\vec{h}\big)\sqrt{\frac{2\lambda_n\big(\vec{h}\big)}{n}}+
\frac{4\lambda_n\big(\vec{h}\big)\cM^2_\infty\big(\vec{h}\big)}{3n},
\end{eqnarray*}
and choosing $\cL=M^2$ and  $\cL_\infty=\cM^2_\infty$, we get from (\ref{eq6-new-prime:proof-th:aux-result1})
\begin{eqnarray}
\label{eq301-new:proof-th:aux-result1}
&&\int_{\bR^d}\bE_f\big\{\chi^p\big(x,\vec{h}\big)\big\}\nu_d(\rd x)
\leq C^{(5)}_pn^{-\frac{p}{2}}e^{(2p+10-b/2)\sum_{j=1}^d(1+\blg_j(\alpha))|\ln(h_j)|},\;\;\;
\forall
\vec{h}\in\cH^d.
\end{eqnarray}

Note that
$
\sigma^{M^2}\big(x,\vec{h}\big)\leq \cM_\infty\big(\vec{h}\big)\sigma\big(x,\vec{h}\big)
$
and, therefore, for any $x\in\bR^d$ and any $\vec{h}\in\cH^d$
\begin{equation}
\label{eq7020:proof-th:aux-result1}
\mU_n\big(x,\vec{h}\big)\leq \cM_\infty\big(\vec{h}\big)U_n\big(x,\vec{h}\big).
\end{equation}
This implies,
\begin{eqnarray*}
\label{eq:proof-th:aux-result2}
\frac{2\lambda_n\big(\vec{h}\big)\widehat{\sigma}^2\big(x,\vec{h}\big)}{n}\leq \frac{2\lambda_n\big(\vec{h}\big)\sigma^2\big(x,\vec{h}\big)}{n}+\frac{2\lambda_n\big(\vec{h}\big)
\cM_\infty\big(\vec{h}\big)U_n\big(x,\vec{h}\big)}{n}
+\frac{2\lambda_n\big(\vec{h}\big)\cM_\infty\big(\vec{h}\big)\chi^*(x,\vec{h}\big)}{n},
\end{eqnarray*}
where we have denoted  $\chi^*(x,\vec{h}\big)=\cM^{-1}_\infty\big(\vec{h}\big)\chi(x,\vec{h}\big)$. Hence
\begin{equation}
\label{eq7021:proof-th:aux-result1}
\widehat{U}_n\big(x,\vec{h}\big)\leq U_n\big(x,\vec{h}\big)+\sqrt{\frac{2\lambda_n\big(\vec{h}\big)
\cM_\infty\big(\vec{h}\big)\big[U_n\big(x,\vec{h}\big)+\chi^*(x,\vec{h}\big)\big]}{n}}.
\end{equation}
By the same reason
\begin{equation}
\label{eq7022:proof-th:aux-result1}
U_n\big(x,\vec{h}\big)\leq \widehat{U}_n\big(x,\vec{h}\big)+\sqrt{\frac{2\lambda_n\big(\vec{h}\big)
\cM_\infty\big(\vec{h}\big)\big[U_n\big(x,\vec{h}\big)+\chi^*\big(x,\vec{h}\big)\big]}{n}}.
\end{equation}
Note that the definition of $\widehat{U}_n\big(x,\vec{h}\big)$ and $U_n\big(x,\vec{h}\big)$ implies that

\begin{equation}
\label{eq7023:proof-th:aux-result1}
\frac{2\lambda_n\big(\vec{h}\big)
\cM_\infty\big(\vec{h}\big)}{n}\leq (3/2)\min\big[\widehat{U}_n\big(x,\vec{h}\big), U_n\big(x,\vec{h}\big)\big].
\end{equation}
Using the inequality $\sqrt{|ab|}\leq 2^{-1}(|ay|+|b/y|)$, $y>0$ we get
 from (\ref{eq7021:proof-th:aux-result1}), (\ref{eq7022:proof-th:aux-result1}) and (\ref{eq7023:proof-th:aux-result1})
\begin{eqnarray*}
\widehat{U}_n\big(x,\vec{h}\big)&\leq& \big(1+\sqrt{3/2}+(3/4)y\big)U_n\big(x,\vec{h}\big)+(2y)^{-1}\chi^*\big(x,\vec{h}\big);
\\
U_n\big(x,\vec{h}\big)&\leq&\big(1+(3/4)y\big)\widehat{U}_n\big(x,\vec{h}\big)+(2y)^{-1}U_n\big(x,\vec{h}\big)+(2y)^{-1}\chi^*\big(x,\vec{h}\big).
\end{eqnarray*}
Choosing $y=1/2$ in the first inequality and $y=1$ in the second we get for any $x\in\bR^d$ and  $\vec{h}\in\cH^d$
\begin{eqnarray}
\label{eq90:proof-th:aux-result2}
\big[\widehat{U}_n\big(x,\vec{h}\big)-3U_n\big(x,\vec{h}\big)\big]_+&\leq& \chi^*(x,\vec{h}\big);
\\
\label{eq91:proof-th:aux-result2}
\big[U_n\big(x,\vec{h}\big)-4\widehat{U}_n\big(x,\vec{h}\big)\big]_+&\leq&\chi^*(x,\vec{h}\big).
\end{eqnarray}
Remembering that  $b=8p+22$ we obtain from (\ref{eq90:proof-th:aux-result2}),  (\ref{eq91:proof-th:aux-result2}), (\ref{eq301-new:proof-th:aux-result1}) and
 (\ref{eq5900:proof-th:aux-result1}) for any $\vec{h}\in\cH^d$
\begin{eqnarray}
\label{eq900:proof-th:aux-result2}
\int_{\bR^d}\bE_f\big[\widehat{U}_n\big(x,\vec{h}\big)-3U_n\big(x,\vec{h}\big)\big]_+^p\nu_d(\rd x)\leq M_\infty^{2p}C^{(5)}_pn^{-\frac{p}{2}}e^{-\sum_{j=1}^d(1+\blg_j(\alpha))|\ln(h_j)|};&&
\\
\label{eq910:proof-th:aux-result2}
\int_{\bR^d}\bE_f\big[U_n\big(x,\vec{h}\big)-4\widehat{U}_n\big(x,\vec{h}\big)\big]_+^p\nu_d(\rd x)\leq M_\infty^{2p}C^{(5)}_pn^{-\frac{p}{2}}e^{-\sum_{j=1}^d(1+\blg_j(\alpha))|\ln(h_j)|}.&&
\end{eqnarray}
The second and third assertions  follow from (\ref{eq900:proof-th:aux-result2}) and (\ref{eq910:proof-th:aux-result2}) with
$
C^\prime_p=M_\infty^{2p}C^{(5)}_p.
$
\epr

\subsection{Proof of Theorem \ref{th:consec-th1-unbounded}}
\label{sec:subsec-proof-concequences}

Let $f\in\bF_{g,\mathbf{u}}(R,D)$.
Introduce the following notations:
$$
a=\Big\{196\big[(c_1\sqrt{c_3})\vee (c_2c_3)\big]\Big\}^{-1},
$$
where $c_1=M_2\sqrt{2D}$,  $c_2=\frac{4M_\infty}{3}$ and
$
c_3=2\max\big\{4\ln(M_\infty), (8p+26)\max_{j=1,\ldots, d}[1+\blg_j(\alpha)]\big\}.
$

\subsubsection{\textsf{Preliminaries}}

Recall that for any locally integrable function $\lambda:\bR^d\to\bR$ its strong maximal function is defined as
\begin{equation}
\label{eq:maximal-function}
 \mM[\lambda](x):= \sup_{H} \frac{1}{\nu_d(H)} \int_H \lambda(t) \rd t,\;\;\;x\in \bR^d,
\end{equation}
where the supremum is taken over all
possible rectangles $H$  in $\bR^d$ with sides parallel to the coordinate
axes, containing point $x$.

It is well known that the strong maximal operator $\lambda\mapsto \mM[\lambda]$
is of the strong $(\mathbf{t},\mathbf{t})$--type for all $1<\mathbf{t}\leq \infty$, i.e.,
if $\lambda\in \bL_\mathbf{t}(\bR^d)$ then $\mM[\lambda]\in \bL_\mathbf{t}(\bR^d)$ and there exists a constant $C_\mathbf{t}$
depending on $\mathbf{t}$ only such that
\begin{equation}
\label{eq:strong-max}
 \big\|\mM[\lambda]\big\|_\mathbf{t} \leq C_\mathbf{t} \|\lambda\|_\mathbf{t},\;\;\;\mathbf{t}\in (1,\infty].
\end{equation}
Let $\mm[\lambda]$ be defined by (\ref{eq:maximal-function}), where, instead of rectangles,  the supremum is taken over all possible cubes $H$  in $\bR^d$ with sides parallel to the coordinate
axes, containing point $x$. Then, it is known that $\lambda\mapsto \mm[\lambda]$ is of the weak $(1,1)$-type, i.e. there exists $C_\mathbf{1}$ depending on $d$ only such that for any $\lambda\in \bL_1(\bR^d)$
\begin{equation}
\label{eq:weak-max-classical}
 \nu_d\Big\{x: \big|\mm[\lambda](x)\big|\geq \mathfrak{z}\Big\} \leq C_\mathbf{1} \mathfrak{z}^{-1}\|\lambda\|_1,\quad\forall\mathfrak{z}>0.
\end{equation}
The results presented below deal with the weak property of the strong maximal function. The following inequality can be found in \cite{Guzman}.
There exists a constant $\mathbf{C}>0$ depending on $d$
only such that
\begin{equation*}
 \nu_d\Big\{x: \big|\mM[\lambda](x)\big|\geq \mathfrak{z}\Big\} \leq \mathbf{C} \int_{\bR^d} \frac{|\lambda(x)|}{\mathfrak{z}}
\bigg\{ 1 +
\bigg(\ln_+\frac{|\lambda(x)|}{\mathfrak{z}}\bigg)^{d-1}\bigg\}\rd x,\quad \mathfrak{z}>0,
\end{equation*}
where for all $z\in\bR$, $\ln_+(z):=\max\{\ln (z),0\}$.

\begin{lemma}
\label{lem:weak-max}
For any given $d\geq 1, R>0$, $ Q>0$ and $\mathbf{q}\in (1,\infty]$
there exists $C(d,\mathbf{q},R,Q)$  such that for any $\lambda\in \bB_{1,d}(R)\cap\bB_{\mathbf{q},d}(Q)$
\begin{equation*}
 \nu_d\Big\{x: \big|\mM[\lambda](x)\big|\geq \mathfrak{z}\Big\} \leq C(d,\mathbf{q},R,Q) \mathfrak{z}^{-1}\big(1+|\ln(\mathfrak{z})|\big)^{d-1},\quad\forall\mathfrak{z}>0.
\end{equation*}
\end{lemma}
The proof of the lemma is an elementary consequence of the aforementioned result and can be omitted.

Recall also the particular case of the Young inequality for weak-type spaces, see \cite{grafakos}, Theorem 1.2.13.
For any $\mathbf{u}\in (1,\infty]$ there exists $C_{\mathbf{u}}>0$ such that for any $\lambda_1\in\bL_1\big(\bR^d\big)$ and $\lambda_2\in\bL_{\mathbf{u},\infty}\big(\bR^d\big)$ one has
\begin{equation}
\label{eq:young-for-weak-spaces}
 \|\lambda_1\star\lambda_2\|_{\mathbf{u},\infty}\leq C_{\mathbf{u}}\|\lambda_1\|_1\|\lambda_2\|_{\mathbf{u},\infty}.
\end{equation}

\paragraph{\textsf{Auxiliary results}} Let us prove several simple facts.
First note that for any $n\geq 3$ for any $\vec{h}\in\cH^d$
\begin{eqnarray}
\label{eq0:proof-add}
&&\lambda_n\big(\vec{h}\big)\leq c_3\Big[\ln{(n)}+\sum_{j=1}^d\big|\ln(\eta_j)\big|\Big].
\end{eqnarray}
Second it is easy to see that for any any $n\geq 3$,
$$
F_n\big(\vec{\eta}\big)\leq F_n\big(\vec{h}\big)\sqrt{l\Big(V_{\vec{\eta}}/V_{\vec{h}}\Big)},\; G_n\big(\vec{\eta}\big)\leq G_n\big(\vec{h}\big)l\Big(V_{\vec{\eta}}/V_{\vec{h}}\Big),\quad \forall  \vec{\eta}, \vec{h}\in (0,\infty)^d:\;\vec{\eta}\geq \vec{h},
$$
 where $l(v)=v^{-1}(1+\ln{v})$. Since $\vec{\eta}\geq \vec{h}$ implies $V_{\vec{\eta}}\geq V_{\vec{h}}$ and $l(v)\le 1$ if $v\ge 1$, we have
\begin{eqnarray}
\label{eq3:proof-add}
F_n\big(\vec{\eta}\big)\leq F_n\big(\vec{h}\big),\quad G_n\big(\vec{\eta}\big)\leq G_n\big(\vec{h}\big),\quad \forall  \vec{\eta}, \vec{h}\in (0,\infty)^d:\;\vec{\eta}\geq \vec{h}.
\end{eqnarray}
Then by (\ref{eq0:proof-add}) and the second inequality in (\ref{eq3:proof-add}), we have:
\begin{eqnarray}
\label{eq5:proof-add}
&&
\sup_{\vec{\eta}\in\cH^d:\;\vec{\eta}\geq \vec{h}}\;\frac{4 M_\infty\lambda_n\big(\vec{\eta}\big)}{3n\prod_{j=1}^d\eta_j(\eta_j\wedge 1)^{\blg_j(\alpha)}}\leq c_2 c_3 G_n\big(\vec{h}\big).
\end{eqnarray}
Now let us establish two bounds for $\|U^*_n\big(\cdot,\vec{h}\big)\|_\infty$.

\smallskip

$\bf 1^0a.\;$ Let $\mathbf{u}=\infty$. We have in view of the second inequality in (\ref{eq:bound-for-supnorm-and-L_2-norm}) for any $\vec{\eta}\in\cH^d$
$$
\sigma\big(x,\vec{\eta}\big)\leq \sqrt{D}\big\|M\big(\cdot,\vec{\eta}\big)\big\|_2\leq M_2\sqrt{D}\prod_{j=1}^d\eta_j^{-\frac{1}{2}}(\eta_j\wedge 1)^{-\blg_j(\alpha)},\quad\forall x\in\bR^d.
$$
It yields for any $x\in\bR^d$ in view of the first inequality in (\ref{eq3:proof-add})
\begin{eqnarray}
\label{eq1:proof-add}
\sup_{\vec{\eta}\in\cH^d:\;\vec{\eta}\geq \vec{h}}\sqrt{\frac{ 2\lambda_n\big(\vec{\eta}\big)\sigma^2\big(x,\vec{\eta}\big)}{n}} \le \sup_{\vec{\eta}\in\cH^d:\;\vec{\eta}\geq \vec{h}} c_1\sqrt{c_3}F_n\big(\vec{\eta}\big)\le  c_1\sqrt{c_3}F_n\big(\vec{h}\big).
\end{eqnarray}
Then gathering (\ref{eq5:proof-add}), (\ref{eq1:proof-add}) and by definition of $a$, we have
\begin{equation}
\label{eq4:proof-add}
\|U^*_n\big(\cdot,\vec{h}\big)\|_\infty\leq (196a)^{-1} \big[F_n\big(\vec{h}\big)
+G_n\big(\vec{h}\big)\big].
\end{equation}


$\bf 1^0b.\;$ Another bound for $\|U^*_n\big(\cdot,\vec{h}\big)\|_\infty$ is available regardless of the value of $\mathbf{u}$. Indeed for any $\vec{\eta}\in\cH^d$ in view of the first inequality in (\ref{eq:bound-for-supnorm-and-L_2-norm})
$$
\sigma\big(x,\vec{\eta}\big)\leq \big\|M\big(\cdot,\vec{\eta}\big)\big\|_\infty\leq M_\infty\prod_{j=1}^d\eta_j^{-1}(\eta_j\wedge 1)^{-\blg_j(\alpha)},\quad\forall x\in\bR^d.
$$
It yields for any $x\in\bR^d$ and any $n\geq 3$
\begin{equation*}
\sup_{\vec{\eta}\in\cH^d:\;\vec{\eta}\geq \vec{h}}\sqrt{\frac{ 2\lambda_n\big(\vec{\eta}\big)\sigma^2\big(x,\vec{\eta}\big)}{n}}
\le \sup_{\vec{\eta}\in\cH^d:\;\vec{\eta}\geq \vec{h}} \sqrt{2c_3} M_\infty\frac{\sqrt{\ln{(n)}+\sum_{j=1}^d\big|\ln(\eta_j)\big|}}{\sqrt{n}\prod_{j=1}^d\eta_j(\eta_j\wedge 1)^{\blg_j(\alpha)}}
\le \sqrt{2c_3 n} M_\infty G_n\big(\vec{h}\big).
\end{equation*}
Then gathering with (\ref{eq5:proof-add}) again, we have
$
\|U^*_n\big(\cdot,\vec{h}\big)\|_\infty\leq \big[(\sqrt{2 c_3 n}M_\infty)\vee (c_2c_3)\big]G_n\big(\vec{h}\big)
$
for any $\vec{h}\in\cH^d$ and, therefore,
\begin{equation}
\label{eq5000:proof-add}
 \inf_{\vec{h}\in\bH}\|U^*_n\big(\cdot,\vec{h}\big)\|_\infty=0.
\end{equation}
To get this it suffices to choose $\vec{h}=(b,\ldots,b)$ and to make $b$ tend to infinity.

\smallskip

$\bf 2^0.\;$  Let now $\mathbf{u}<\infty$. Let us prove that for any $\mz>0$, $\mathbf{s}\in\{1,\mathbf{u}\}$ and any $f\in\bF_{g,\mathbf{u}}(R,D)$
\begin{eqnarray}
\label{eq6:proof-add}
\nu_d\bigg(x\in\bR^d:\; \sup_{\vec{\eta}\in\cH^d:\;\vec{\eta}\geq \vec{h}}\;\cU_n\big(x,\vec{\eta},f\big)\geq \mz\bigg)\leq c_5\big[\widetilde{D}\mz^{-2}F_n^2\big(\vec{h}\big)\big]^{\mathbf{s}},
\end{eqnarray}
where we have put $\cU^2_n\big(\cdot,\vec{\eta},f\big)=2n^{-1}\lambda_n\big(\vec{\eta}\big)\sigma^2\big(\cdot,\vec{\eta}\big)$ and $\widetilde{D}=1$ if $\mathbf{s}=1$
and $\widetilde{D}=D$ if $\mathbf{s}=\mathbf{u}$.

Indeed, if $\mathbf{s}=1$, applying the Markov inequality,  we obtain in view of the second inequality in (\ref{eq:bound-for-supnorm-and-L_2-norm}) for any $\vec{\eta}\in\cH^d$
\begin{eqnarray}
\label{eq7:proof-add}
&&\nu_d\Big(x\in\bR^d:\; \cU_n\big(x,\vec{\eta},f\big)\geq \mz\Big)\leq 2(n\mz^{2})^{-1}\lambda_n\big(\vec{\eta}\big)\int_{\bR^d}\sigma^2\big(x,\vec{\eta}\big)\nu_d(\rd x)
\\&&
= 2(n\mz^{2})^{-1}\lambda_n\big(\vec{\eta}\big)\big\|M\big(\cdot,\vec{\eta}\big)\big\|^2_2
\leq
2M_2^2(n\mz^{2})^{-1}\frac{\lambda_n\big(\vec{\eta}\big)}{\prod_{j=1}^d\eta_j(\eta_j\wedge 1)^{2\blg_j(\alpha)}}
\leq c_6\mz^{-2}F_n^2\big(\vec{\eta}\big).
\nonumber
\end{eqnarray}
Here we have put $c_6=2M_2^2c_1^2 c_3$ and to get the last inequality we have used (\ref{eq0:proof-add}).

To get the similar result  if $\mathbf{s}=\mathbf{u}$ we remark that $\sigma^2\big(\cdot,\vec{\eta}\big)=M^2\big(\cdot,\vec{\eta}\big)\star\mathfrak{p}(\cdot)$ and that $M^2\big(\cdot,\vec{\eta}\big)\in \bL_1\big(\bR^d\big)$ in view of the second inequality in (\ref{eq:bound-for-supnorm-and-L_2-norm}). It remains to note
that $f\in\bF_{g,\mathbf{u}}(R,D)$ implies $\mathfrak{p}\in\bB_{\mathbf{u},d}^{(\infty)}(D)$ and to apply the inequality (\ref{eq:young-for-weak-spaces}).

It yields together with the second inequality in (\ref{eq:bound-for-supnorm-and-L_2-norm}) for any $\vec{\eta}\in\cH^d$
\begin{eqnarray}
\label{eq7000:proof-add}
&&\nu_d\Big(x\in\bR^d:\; \cU_n\big(x,\vec{\eta},f\big)\geq \mz\Big)\leq
\big[c_6C_{\mathbf{u}}D\mz^{-2}F_n^2\big(\vec{\eta}\big)\big]^{\mathbf{u}}.
\end{eqnarray}
Thus, denoting $\widetilde{C}=1$ if $\mathbf{s}=1$
and $\widetilde{C}=C_{\mathbf{u}}$ if $\mathbf{s}=\mathbf{u}$, we get from (\ref{eq7:proof-add}) and (\ref{eq7000:proof-add})
\begin{eqnarray}
\label{eq8:proof-add}
&&\nu_d\bigg(x\in\bR^d:\; \sup_{\vec{\eta}\in\cH^d:\;\vec{\eta}\geq \vec{h}}\;\cU_n\big(x,\vec{\eta},f\big)\geq \mz\bigg)
\leq \big[c_6\widetilde{C}\widetilde{D}\mz^{-2}\big]^{\mathbf{s}}\sum_{\vec{\eta}\in\cH^d:\;\vec{\eta}\geq \vec{h}}\;F_n^{2\mathbf{s}}\big(\vec{\eta}\big).
\end{eqnarray}
It remains to note that since $\vec{\eta},\vec{h}\in\cH^d$ and  $\vec{\eta}\geq\vec{h}$ we can write $\eta_j=e^{m_j}h_j$ with $m_j\geq 0$ for any $j=1,\ldots, d$. It yields together with the first inequality in (\ref{eq3:proof-add})
$$
\sum_{\vec{\eta}\in\cH^d:\;\vec{\eta}\geq \vec{h}}\;F_n^{2\mathbf{s}}\big(\vec{\eta}\big)\leq F_n^{2\mathbf{s}}\big(\vec{h}\big)\sum_{(m_1,\ldots, m_d)\in\bN^d}\Big(1+\sum_{j=1}^dm_j\Big)^{\mathbf{s}}e^{-\mathbf{s}\sum_{j=1}^dm_j}=:c_7F_n^{2\mathbf{s}}\big(\vec{h}\big).
$$
Hence, (\ref{eq6:proof-add}) with $c_5=c_7[c_6\widetilde{C}]^{\mathbf{s}}$ follows from (\ref{eq8:proof-add}).

\smallskip

$\bf 3^0.\;$ Let $c_{K}\geq 1$ be such that $\text{supp}(K)\subset [-c_K,c_K]^d$. We have
$$
\big|S_{\vec{h}}(x,f)\big|=\Big|\int_{\bR^d}K_{\vec{h}}(t-x)f(t)\nu_d(\rd t)\Big|\leq (2c_K)^d\|K\|^d_\infty \mM[|f|](x),\quad\forall \vec{h}\in (0,\infty)^d.
$$
If $\vec{h}=(h,\ldots,h), h\in(0,\infty)$, the latter inequality  holds with $\mm[|f|]$ instead of $\mM[|f|]$.
Thus,
\begin{eqnarray}
\label{eq10:proof-add}
&&\sup_{\vec{h}\in\bH}\big|\cB_{\vec{h}}(x,f)\big|\leq 3(2c_K)^d\|K\|_\infty \mM_{\bH}[|f|](x)+|f(x)|,\quad\forall  x\in\bR^d,
\end{eqnarray}
where we have denoted $\mM_{\bH}=\mM$ if $\bH=\cH^d$ and $\mM_{\bH}=\mm$ if $\bH=\cH^d_{\text{isotr}}$.

Moreover, we deduce from (\ref{eq10:proof-add}) and (\ref{eq5000:proof-add}) putting $T_{\vec{h}}(x,f)=\cB_{\vec{h}}(x,f)+49U^*_n\big(\cdot,\vec{h}\big)$ that
\begin{eqnarray}
\label{eq1000:proof-add}
&&\inf_{\vec{h}\in\bH}\big|T_{\vec{h}}(x,f)\big|\leq 3(2c_K)^d\|K\|_\infty \mM_{\bH}[|f|](x)+|f(x)|.
\end{eqnarray}


\subsubsection{\textsf{Proof of the theorem}}

For any $v>0$ set $\cC_v(f)=\big\{x\in\bR^d:\;\mathbf{T}(x,f)\geq v\big\}$, where we have put $\mathbf{T}(x,f)=\inf_{\vec{h}\in\bH}|T_{\vec{h}}(x,f)|$.
For any given $\overline{\Bv}>0$ one obviously has
\begin{eqnarray}
\label{eq14:proof-add}
\|\mathbf{T}(\cdot,f)\|^p_p\leq p\int_{0}^{\overline{\Bv}}v^{p-1}\nu_d\big(\cC_{v}(f)\big)\rd v+
\int_{\cC_{\overline{\Bv}}(f)}|\mathbf{T}(x,f)|^p\nu_d(\rd x)
\end{eqnarray}
Denoting $\cW_v(\vec{h},f)=\{x\in\bR^d:\; 49U^*_n\big(x,\vec{h}\big) \geq 2^{-1}v\}$ we obviously have for any $\vec{h}\in\bH$ and $v>0$
\begin{eqnarray}
\label{eq15:proof-add}
\nu_d\big(\cC_v(f)\big)&\leq& \nu_d\big(\cA(\vec{h},f,v)\big)+\nu_d\big(\cW_v(\vec{h},f)\big);
\\*[2mm]
\label{eq1500:proof-add}
|\mathbf{T}(x,f)|^p\mathrm{1}_{\cC_v(f)}(x)&\leq& 2^p\big|\cB_{\vec{h}}(x,f)\big|^p\mathrm{1}_{\cA(\vec{h},f,v)}
+98^p\big|U^*_n\big(x,\vec{h}\big)\big|^p\mathrm{1}_{\cW_{v}(\vec{h},f)}(x);
\\*[2mm]
\label{eq15000:proof-add}
\nu_d\big(\cC_v(f)\big)&\leq& \nu_d\Big(x\in\bR^d:\; 3(2c_K)^d\|K\|_\infty \mM_{\bH}[|f|](x)+|f(x)|>v \Big).
\end{eqnarray}
The last inequality follows from (\ref{eq1000:proof-add}). Set $\cU^*_n\big(x,\vec{h},f\big)=\sup_{\vec{\eta}\in\cH^d:\;\vec{\eta}\geq \vec{h}}\;\cU_n\big(x,\vec{\eta},f\big)$.

\smallskip

$\bf 1^0.\;$
Noting that  $U^*_n\big(x,\vec{h}\big)\leq \cU^*_n\big(x,\vec{h},f\big)+(196a)^{-1}G_n\big(\vec{h}\big)$ in view of (\ref{eq5:proof-add}), we get
\begin{eqnarray}
\label{eq1001:proof-add}
&& \cW_v(\vec{h},f)\subseteq\Big\{x\in\bR^d:\; 49\cU^*_n\big(x,\vec{h}\big)\geq 4^{-1}v\Big\}:=\widetilde{\cW}_v(\vec{h},f),\quad\forall \vec{h}\in\mH(v).
\end{eqnarray}
Applying (\ref{eq6:proof-add}) with $\mathbf{s}=1$ we deduce from (\ref{eq15:proof-add}) that
\begin{eqnarray*}
&&\nu_d\big(\cC_v(f)\big)\leq \nu_d\big(\cA(\vec{h},f,v)\big)+196^{2}c_5 v^{-2} F_n^2\big(\vec{h}\big), \quad\forall \vec{h}\in\mH(v).
\end{eqnarray*}
Noting that the left hand side of the latter inequality is independent of $\vec{h}$ we get
\begin{eqnarray}
\label{eq151:proof-add}
&&\nu_d\big(\cC_v(f)\big)\leq \max[1, 196^{2}c_5]\Lambda(v,f).
\end{eqnarray}

\smallskip

$\bf 2^0.\;$ Let us establish the following bounds, where $c_9$ is given in the paragraph $\bf 2^0b.$ below.

\noindent For any $\mathbf{u}\in[1,\infty]$,
\begin{eqnarray}
\label{eq18:proof-add}
&&\nu_d\big(\cC_v(f)\big)\leq \max[1,c_5196^{2}, c_5196^{2\mathbf{u}}D^{\mathbf{u}}a^{2\mathbf{u}}]\{\Lambda(v,f)\wedge\Lambda(v,f,\mathbf{u})\},\quad \forall v>0.
\end{eqnarray}
and for any $\mathbf{u}\in(p/2,\infty]$,
\begin{eqnarray}
\label{eq19:proof-add}
\int_{\cC_{v}(f)}|\mathbf{T}(x,f)|^p\nu_d(\rd x)\leq \max[2^p,98^pc_9]\Lambda_{p}\big(v,f,\mathbf{u}\big), \quad \forall v>0.
\end{eqnarray}

\quad $\bf 2^0a.\;$ Let $\mathbf{u}=\infty$. We remark that  $ \cW_v(\vec{h},f)=\emptyset$ for any $\vec{h}\in\mH(v,2)$ in view of (\ref{eq4:proof-add}).
Thus, we deduce from (\ref{eq15:proof-add}),  (\ref{eq1500:proof-add}) and (\ref{eq:case-mathbf(u)=infty}), taking into account that the left hand sides of both inequalities are independent of $\vec{h}$
\begin{eqnarray}
\label{eq152:proof-add}
&&\nu_d\big(\cC_v(f)\big)\leq \Lambda(v,f,\infty),\qquad \int_{\cC_{v}(f)}|\mathbf{T}(x,f)|^p\nu_d(\rd x)\leq \Lambda_{p}(v,f,\infty).
\end{eqnarray}
This inequality and (\ref{eq151:proof-add}) ensure that (\ref{eq18:proof-add}) and (\ref{eq19:proof-add}) hold if $\mathbf{u}=\infty$.

\smallskip

\quad $\bf 2^0b.\;$ Let $\mathbf{u}<\infty$. Applying (\ref{eq6:proof-add}) with $\mathbf{s}=\mathbf{u}$, we obtain in view of (\ref{eq1001:proof-add})
\begin{eqnarray*}
&&\nu_d\big(\cW_v(\vec{h},f)\big)\leq c_5196^{2\mathbf{u}}D^{\mathbf{u}}v^{-2\mathbf{u}}F_n^{2\mathbf{u}}\big(\vec{h}\big)\leq c_5196^{2\mathbf{u}}D^{\mathbf{u}}a^{2\mathbf{u}}z^{-\mathbf{u}},\quad \forall \vec{h}\in\mH(v,z)
\end{eqnarray*}
It yields together with (\ref{eq15:proof-add})
\begin{eqnarray}
\label{eq16:proof-add}
&&\nu_d\big(\cC_v(f)\big)\leq \max[1, c_5196^{2\mathbf{u}}D^{\mathbf{u}}a^{2\mathbf{u}}]\Lambda(v,f,\mathbf{u}).
\end{eqnarray}
This inequality and (\ref{eq151:proof-add}) ensure that (\ref{eq18:proof-add}) holds if $\mathbf{u}<\infty$.


\smallskip

What is more, we have in view of (\ref{eq5:proof-add}) and (\ref{eq1001:proof-add}) for any $\vec{h}\in\mH(v)$
$$
\big|U^*_n\big(x,\vec{h}\big)\big|^p\mathrm{1}_{\cW_{v}(\vec{h},f)}\leq 2^{p}\big|\cU^*_n\big(x,\vec{h},f\big)\big|^p\mathrm{1}_{\widetilde{\cW}_{v}(\vec{h},f)}
$$
Moreover, applying (\ref{eq6:proof-add}) with $\mathbf{s}=\mathbf{u}$, we have for any $y>0$ and  $\vec{h}\in\mH(v,z)$
$$
\nu_d\big(\widetilde{\cW}_{y}(\vec{h},f)\big)\leq c_5196^{2\mathbf{u}}D^{\mathbf{u}}y^{-2\mathbf{u}}F_n^{2\mathbf{u}}\big(\vec{h}\big)\leq
c_5196^{2\mathbf{u}}D^{\mathbf{u}}y^{-2\mathbf{u}}(av)^{2\mathbf{u}}z^{-\mathbf{u}}.
$$
Hence, if additionally $\mathbf{u}>p/2$, we have for any $\vec{h}\in\mH(v,z)$
\begin{eqnarray*}
\int_{\cW_{v}(\vec{h},f)}\big|U^*_n\big(x,\vec{h}\big)\big|^p\nu_d(\rd x)&\le& 2^pp\int_{v}^\infty y^{p-1}\nu_d\big(\widetilde{\cW}_{y}(\vec{h},f)\big)\rd y
\\
&=& c_5196^{2\mathbf{u}}D^{\mathbf{u}}a^{2\mathbf{u}}2^{p}pv^{2\mathbf{u}}z^{-\mathbf{u}}\int_{v}^\infty y^{p-1-2\mathbf{u}}\rd y
=:c_9v^{p}z^{-\mathbf{u}}.
\end{eqnarray*}
This yields together with (\ref{eq1500:proof-add})
\begin{eqnarray}
\label{eq17:proof-add}
\int_{\cC_{v}(f)}|\mathbf{T}(x,f)|^p\nu_d(\rd x)\leq \max[2^p,98^pc_9]\Lambda_{p}\big(v,f,\mathbf{u}\big).
\end{eqnarray}
This inequality ensures that (\ref{eq19:proof-add}) holds if $\mathbf{u}<\infty$.

\smallskip

$\bf 3^0.\;$ Recall that $f\in\bF_g(R)$ implies that $f\in\bB_{\mathbf{1},d}(R)$. Since additionally  $f\in\bB_{\mathbf{q},d}(D)$, $\mathbf{q}>1$, Lemma \ref{lem:weak-max} as well as (\ref{eq:weak-max-classical}) is applicable and we obtain in view of  (\ref{eq15000:proof-add})
$$
\nu_d\big(\cC_v(f)\big)\leq c_{10}v^{-1}(1+|\ln{v}|)^{t(\bH)},\quad \forall v>0.
$$
It yields for any $\overline{\Bv}>0$ and $p>1$
\begin{eqnarray}
\label{eq190:proof-add}
&&
p\int_{0}^{\underline{\Bv}}v^{p-1}\nu_d\big(\cC_{v}(f)\big)\rd v\leq c_{10} p\int_{0}^{\underline{\Bv}}v^{p-2}(1+|\ln{v}|)^{t(\bH)}\rd v\leq c_{11}\underline{\Bv}^{p-1}(1+|\ln{\underline{\Bv}}|)^{t(\bH)}.
\end{eqnarray}
In the case of $t(\bH)=0$ the last inequality is obvious and if $t(\bH)=d-1$ it follows by integration by parts.
The assertion of the theorem follows now from  (\ref{eq14:proof-add}), where the bound  (\ref{eq190:proof-add}) is used for any $v<\underline{\Bv}$,
the estimate (\ref{eq18:proof-add}) for any $v\in[\underline{\Bv}, \overline{\Bv}]$ and the bound (\ref{eq19:proof-add}) with $v=\overline{\Bv}$.
\epr

\subsection{Proof of Theorem \ref{th:minimax-abstract}}
\label{sec:subsec-proof-th:minimax-abstract}

The proof of the theorem is based essentially on some auxiliary statements formulated in Section \ref{sec:subsubsec-Auxiliary results-proof-th:minimax-abstract} below.

Some properties related to the kernel approximation of the underlying function $f$ are summarized in Lemma \ref{lem:poinwise-approx-bound} and in formulae (\ref{eq:Young-restricted}). The results presented in Lemma \ref{lem:weak-max} and in formulae (\ref{eq:strong-max-partial-operator}) deal with the properties of the strong maximal function.
In the subsequent proof
$c_1,c_2, \ldots$,
stand for constants depending only on $\vec{s}, \vec{q},  g, \cK,d$, $R,D, \mathbf{u}$ and $\mathbf{q}$.

\subsubsection{\textsf{Auxiliary results}}
\label{sec:subsubsec-Auxiliary results-proof-th:minimax-abstract}

Let $\mJ$ denote the set of all the subsets of $\{1,\ldots d\}$ endowed with the empty set $\emptyset$. For any $J\in\mJ$ and  $y\in\bR^d$ set $y_J=\{y_j,\;j\in J\}\in\bR^{|J|}$ and we will write $y=\big(y_J,y_{\bar{J}}\big)$, where as usual $\bar{J}=\{1,\ldots d\}\setminus J$.

For any
$j=1,\ldots,d$ introduce the $d\times d$ matrix $\mathbf{E}_j=(\mathbf{0},\ldots,\mathbf{e}_j,\ldots,\mathbf{0})$ where, recall, $(\mathbf{e}_1,\dots,\mathbf{e}_d)$ denotes the canonical basis of $\mathbb{R}^d$. Set also $\mathbf{E}[J]=\sum_{j\in J}\mathbf{E}_j$. Later on $\mathbf{E}_0=\mathbf{E}[\emptyset]$ denotes the matrix with zero entries.

To any  $J\in\mJ$  and  any $\lambda:\bR^d\to\bR$   associate the function
$$
\lambda_J\big(y_J,z_{\bar{J}}\big)=\lambda\big(z+\mathbf{E}[J](y-z)\big),\quad y,z\in\bR^d,
$$
with the obvious agreement $\lambda_{J}\equiv\lambda$ if $J=\{1,\ldots d\}$, which is always the case if $d=1$.

For any  $\vec{h}\in\cH^d$ and $J\subseteq\{1,\ldots d\}$ set
$
K_{\vec{h},J}(u_J)=\prod_{j\in J}h^{-1}_j\cK\big(u_j/h_j\big)
$
and define for any $y\in\bR^d$
$$
\big[K_{\vec{h}}\circ\lambda\big]_{J}(y)=\int_{\bR^{|\bar{J}|}}K_{\vec{h},\bar{J}}(u_{\bar{J}}-y_{\bar{J}})\lambda\big(y_J,u_{\bar{J}}\big)
\nu_{|\bar{J}|}\big(\rd u_{\bar{J}}\big),
$$ where $\nu_{|\bar{J}|}$ is the Lebesgue measure on $\bR^{|\bar{J}|}$.
For any $\vec{h},\vec{\eta}\in\cH^d$ set
$$
B_{\vec{h},\vec{\eta}}(x,f)=|S_{\vec{h}\vee \vec{\eta}}(x,f)-S_{ \vec{\eta}}(x,f)|.
$$

\begin{lemma}
\label{lem:poinwise-approx-bound}
Let Assumption \ref{ass2:on-kernel-deconvolution-part1} hold. One can find $k \in \{1,\ldots d\}$ and a collection of indexes
 $\big\{j_1<j_2<\cdots<j_k\big\}\in \{1,\ldots,d\}$ such that for any  $x\in\bR^d$ and any  $f:\bR^d\to\bR$
\begin{eqnarray*}
B_{\vec{h},\vec{\eta}}(x,f)&\leq& \sum_{l=1}^k\Big(\left[\big|K_{\vec{h}\vee\vec{\eta}}\big|\circ b_{h_{j_l},f,j_l} \right]_{J_l}(x)+
\left[\big|K_{\vec{\eta}}\big|\circ b_{h_{j_l},f,j_l}\right]_{J_l}(x)\Big);
\\
B_{\vec{h}}(x,f)&\leq& \sum_{l=1}^k
\left[\big|K_{\vec{h}}\big|\circ b_{h_{j_l},f,j_l} \right]_{J_l}(x),\qquad J_l=\{j_1,\ldots, j_l\}.
\end{eqnarray*}

\end{lemma}

The proof of the lemma can be found in \cite{lepski15},  Lemma 2.

Also, let us mention the following bound which is a trivial consequence of the Young inequality and the Fubini theorem.
If $\lambda\in\bL_\mathbf{t}(\bR^d)$ then for any $\mathbf{t}\in[1,\infty]$
\begin{equation}
\label{eq:Young-restricted}
\sup_{J\in\mJ}\big\|\big[K_{\vec{h}}\circ\lambda\big]_{J}\big\|_{\mathbf{t}}\leq \|\cK\|^{d}_{1}\|\lambda\|_{\mathbf{t}},\quad\forall \vec{h}\in\cH^d.
\end{equation}

To any $J\in\mJ$ and  any locally integrable function $\lambda:\bR^d\to\bR_+$ we associate the operator
\begin{equation}
\label{eq:strong-max-partial-operator}
\mM_{J}[\lambda](x)=\sup_{H_{|\bar{J}|}}\frac{1}{\nu_{|\bar{J}|}(H_{|\bar{J}|})}\int_{H_{|\bar{J}|}}\lambda\big(t+\mathbf{E}[J][x-t]\big)
\nu_{|\bar{J}|}(\rd t_{\bar{J}})
\end{equation}
where
 the supremum is taken over all hyper-rectangles in $\bR^{|\bar{J}|}$ containing $x_{\bar{J}}=(x_j, j\in \bar{J})$ and with sides parallel to the axis.

As we see $\mM_J[\lambda]$ is the strong maximal operator applied to the function obtained from $\lambda$ by fixing the coordinates whose indices belong to $J$. It is obvious  that $\mM_{\emptyset}[\lambda]\equiv \mM[\lambda]$ and $\mM_{\{1,\ldots,d\}}[\lambda]\equiv\lambda$.

The following result is a direct consequence of (\ref{eq:strong-max}) and of the Fubini theorem.
For any  $\mathbf{t}\in(1,\infty]$  there exists $\mathbf{C}_{\mathbf{t}}$ such that for any  $\lambda\in\bL_\mathbf{t}\big(\bR^d)$
\begin{equation}
\label{eq:strong-max-partial}
\sup_{J\in\mJ}\big\|\mM_J[\lambda]\big\|_{\mathbf{t}}\leq \mathbf{C}_{\mathbf{t}}\|\lambda\|_{\mathbf{t}}.
\end{equation}
Obviously this inequality holds if $\mathbf{t}=\infty$ with $\mathbf{C}_{\infty}=1$.

\subsubsection{Proof of the theorem}

$\mathbf{1^0.}\;$ We start with the following obvious observation.
For any $\lambda:\bR^d\to\bR_+$, $\vec{u}\in\bR^d$ and $J\in\mJ$
\begin{equation}
\label{eq1:proof-prop:measure-of-bias-deconv}
[K_{\vec{u}}\circ\lambda]_J(x)\leq (2c_\cK\|\cK\|_\infty)^{d}\mM_J[\lambda](x),\quad\forall x\in\bR^d.
\end{equation}
Putting $C_1=(2c_\cK\|\cK\|_\infty)^{d}$
we get  for any $\vec{h},\vec{\eta}\in\cH^d$ and $x\in\bR^d$ in view of (\ref{eq1:proof-prop:measure-of-bias-deconv}) and   assertions of Lemma \ref{lem:poinwise-approx-bound} that
\begin{eqnarray*}
B_{\vec{h},\vec{\eta}}(x,f)\leq 2C_1\sum_{j=1}^d\sup_{J\in\mJ}\mM_{J}\big[b_{h_{j},f,j}\big](x),\quad B_{\vec{h}}(x,f)\leq C_1\sum_{j=1}^d\sup_{J\in\mJ}\mM_{J}\big[b_{h_{j},f,j}\big](x).
\end{eqnarray*}
Thus noting that the right hand side of the first inequality above is independent of $\vec{\eta}$, we obtain
\begin{eqnarray}
\label{eq2001:proof-prop:measure-of-bias-deconv}
\cB_{\vec{h}}(x,f)\leq 5C_1\sum_{j=1}^d\sup_{J\in\mJ}\mM_{J}\big[b_{h_{j},f,j}\big](x),\quad\forall x\in\bR^d,\; \forall\vec{h}\in\cH^d.
\end{eqnarray}
Applying  (\ref{eq:strong-max-partial}) with $\mathbf{t}=\infty$,
we have for any $v>0$ in view of  the definition of  $J(\vec{h},v)$
\begin{eqnarray}
\label{eq2000:proof-prop:measure-of-bias-deconv}
\cB_{\vec{h}}(x,f)&\leq& 5C_1\sum_{j\in \bar{J}(\vec{h},v)}\sup_{J\in\mJ}\mM_{J}\big[b_{h_{j},f,j}\big](x)+5C_1\sum_{j\in J(\vec{h},v)}
\sup_{J\in\mJ}\Big\|\mM_{J}\big[b_{h_{j},f,j}\big]\Big\|_\infty
\nonumber\\
&\leq&  5C_1\sum_{j\in\bar{J}(\vec{h},v)}\sup_{J\in\mJ}\mM_{J}\big[b_{h_{j},f,j}\big](x)+5C_1\sum_{j\in J(\vec{h},v)}\mathbf{B}_{j,\infty,\bF}\big(h_{j}\big)
\nonumber\\
&\leq& 5C_1\sum_{j\in\bar{J}(\vec{h},v)}\sup_{J\in\mJ}\mM_{J}\big[b_{h_{j},f,j}\big](x)+4^{-1}v, \quad \forall f\in\bF.
\end{eqnarray}
We obtain for any $f\in\bF$, $v>0$ and  $\vec{s}=(s_1,\ldots,s_d)\in (1,\infty)^d$, applying consecutively  the Markov inequality and    (\ref{eq:strong-max-partial}) with $\mathbf{t}=s_j$,
\begin{eqnarray}
\label{eq2:proof-prop:measure-of-bias-deconv}
\nu_d\left\{\cA\big(\vec{h},f,v\big)\right\}&\leq&\nu_d \Big(\cup_{J\in\mJ}\cup_{j\in\bar{J}(\vec{h},v)}\Big\{x: 5C_1\mM_{J}\big[b_{h_{j},f,j}\big](x)\geq (4d)^{-1}v\Big\}\Big)
\nonumber\\
&\leq& c_1\sum_{j\in\bar{J}(\vec{h},v)} v^{-s_j}\big\|b_{h_{j},f,j}\big\|_{s_j}^{s_j}
\le c_1\sum_{j\in\bar{J}(\vec{h},v)} v^{-s_j}\Big[\mathbf{B}_{j,s_j,\bF}\big(h_{j}\big)\Big]^{s_j}.
\end{eqnarray}
Noting that the right hand side of the latter inequality is independent of $f$ and the left hand side is independent of $\vec{s}$, we get
\begin{eqnarray}
\label{eq222:proof-prop:measure-of-bias-deconv}
&&c_1^{-1}\sup_{f\in\bF}\{\Lambda(v,f)\wedge\Lambda(v,f,\mathbf{u})\}\leq \blL_{\vec{s}}(v,\bF,\mathbf{u})
\wedge\blL_{\vec{s}}(v,\bF),  \quad \forall v>0,\;\;\vec{s}\in (1,\infty)^d.
\end{eqnarray}

$\mathbf{2^0.}\;$ Note also that in view of (\ref{eq2000:proof-prop:measure-of-bias-deconv}), we have for any  $v>0$
\begin{eqnarray}
\label{eq888:proof-add}
&&\int_{\cA\big(\vec{h},f,v\big)} \big|\cB_{\vec{h}}(x,f)\big|^{p}\nu_d(\rd x)
\nonumber\\&&
\leq c_2\int_{\cA\big(\vec{h},f,v\big)} \bigg|\sum_{j\in\bar{J}(\vec{h},v)}\sup_{J\in\mJ}\mM_{J}\big[b_{h_{j},f,j}\big](x)\bigg|^{p}\nu_d(\rd x)+c_3v^{p}\nu_d\left\{\cA\big(\vec{h},f,v\big)\right\}
\nonumber\\
&&\leq c_4\Bigg[\sum_{j\in\bar{J}(\vec{h},v)}\int_{\cA\big(\vec{h},f,v\big)}\Big|\sup_{J\in\mJ}\mM_{J}\big[b_{h_{j},f,j}\big](x)\Big|^{p}\nu_d(\rd x)+ v^p\nu_d\left\{\cA\big(\vec{h},f,v\big)\right\}\Bigg].
\end{eqnarray}
For any $v>0$ and $j=1,\ldots,d,$ introduce
$$
\mathfrak{A}_j(v)= \Big\{x\in\bR^d:\; \sup_{J\in\mJ}\mM_{J}\big[b_{h_{j},f,j}\big](x)\geq (40C_1)^{-1}v\Big\},\quad \cA_j(v)=\cA\big(\vec{h},f,v\big)\cap\bar{\mA}_j(v).
$$
Noting that in view of (\ref{eq2000:proof-prop:measure-of-bias-deconv}) for any $v>0$ and any $j\in\bar{J}(\vec{h},v)$
\begin{eqnarray*}
\cA_j(v)&\subseteq& \Bigg\{x\in\bR^d:\; 5C_1\sum_{k\in\bar{J}(\vec{h},v),\;k\neq j}\sup_{J\in\mJ}\mM_{J}\big[b_{h_{j},f,k}\big](x)\geq v/8\Bigg\}
\\
&\subseteq& \Bigg\{x\in\bR^d:\; 5C_1\sum_{k\in\bar{J}(\vec{h},v)}\sup_{J\in\mJ}\mM_{J}\big[b_{h_{j},f,k}\big](x)\geq v/8\Bigg\}=:\cA^*\big(\vec{h},f,v\big),
\end{eqnarray*}
we deduce from (\ref{eq888:proof-add}) that for any $\vec{q}\in [p,\infty)^d$
\begin{eqnarray}
\label{eq8888:proof-add}
&&\int_{\cA\big(\vec{h},f,v\big)} \big|\cB_{\vec{h}}(x,f)\big|^{p}\nu_d(\rd x)\leq c_4\sum_{j\in\bar{J}(\vec{h},v)}\int_{\mA_j(v)}\Big|\sup_{J\in\mJ}\mM_{J}\big[b_{h_{j},f,j}\big](x)\Big|^{p}\nu_d(\rd x)
\nonumber\\
&&\qquad\qquad\qquad\qquad\qquad\qquad\qquad+c_5v^p\Big[\nu_d\left\{\cA^*\big(\vec{h},f,v\big)\right\} + \nu_d\left\{\cA\big(\vec{h},f,v\big)\right\}\Big]
\nonumber\\
&&\leq c_6\sum_{j\in\bar{J}(\vec{h},v)}v^{p-q_j}\Big\|\sup_{J\in\mJ}\mM_{J}\big[b_{h_{j},f,j}\big]\Big\|^{q_j}_{q_j}
+c_5v^p\Big[\nu_d\left\{\cA^*\big(\vec{h},f,v\big)\right\} + \nu_d\left\{\cA\big(\vec{h},f,v\big)\right\}\Big].
\end{eqnarray}
It remains to note that similarly (\ref{eq2:proof-prop:measure-of-bias-deconv}) for any $\vec{s}\in (1,\infty)^d$
$$
\nu_d\left\{\cA^*\big(\vec{h},f,v\big)\right\}\leq c_7 \sum_{j\in\bar{J}(\vec{h},v)} v^{-s_j}\big\|b_{h_{j},f,j}\big\|_{s_j}^{s_j}
$$
and to apply
(\ref{eq:strong-max-partial}) with $\mathbf{t}=q_j$  to the each term in the sum appeared in (\ref{eq8888:proof-add}). All of this together with (\ref{eq2:proof-prop:measure-of-bias-deconv}), applied with $\vec{s}=\vec{q}$ yields for any  $v>0$ and  $\vec{q}\in [p,\infty)^d$
\begin{eqnarray*}
\label{eq30:proof-prop:measure-of-bias-deconv}
\int_{\cA\big(\vec{h},f,v\big)} \big|\cB_{\vec{h}}(x,f)\big|^{p}\nu_d(\rd x)&\leq& c_9\sum_{j\in\bar{J}(\vec{h},v)} v^{p-q_j}\big\|b_{h_{j},f,j}\big\|_{q_j}^{q_j}\le c_9\sum_{j\in\bar{J}(\vec{h},v)} v^{p-q_j}\Big[\mathbf{B}_{j,q_j,\bF}\big(h_{j}\big)\Big]^{q_j}.
\end{eqnarray*}
Noting that the right hand side of the latter inequality is independent of $f$ and the left hand side is independent of $\vec{q}$, the we get
\begin{eqnarray}
\label{eq22222:proof-prop:measure-of-bias-deconv}
&&\sup_{f\in\bF}\Lambda_p(v,f,\mathbf{u})\leq c_9v^p\blL_{\vec{q}}(v,\bF,\mathbf{u}),  \quad \forall v>0,\;\;\vec{q}\in [p,\infty)^d.
\end{eqnarray}
The first assertion of the theorem follows from (\ref{eq222:proof-prop:measure-of-bias-deconv}), (\ref{eq22222:proof-prop:measure-of-bias-deconv}) and Theorem
\ref{th:consec-th1-unbounded}.

\smallskip

$\mathbf{3^0.}\;$
Remark that in view of (\ref{eq10:proof-add}) and (\ref{eq:strong-max}) $f\in\bB_{\mathbf{q},d}(D)$ implies
\begin{eqnarray}
\label{eq9:proof-add}
&&\big\|\cB_{\vec{h}}(\cdot,f)\big\|_\mathbf{q}\leq \big[3(2c_\cK)^d\|\cK\|^d_\infty C_{\mathbf{q}}+1\big]D,\quad\forall \vec{h}\in (0,\infty)^d,
\end{eqnarray}
where $C_{\mathbf{q}}$ is the constant which appeared in (\ref{eq:strong-max}).
Hence  for any  $v>0$ and $\mathbf{q}\in [p,\infty)$
\begin{eqnarray}
\label{eq251:proof-add}
&&\int_{\cA\big(\vec{h},f,v\big)} \big|\cB_{\vec{h}}(x,f)\big|^{p}\nu_d(\rd x)\leq 2^{\mathbf{q}-p} v^{p-\mathbf{q}}\big\|\cB_{\vec{h}}(\cdot,f)\big\|_\mathbf{q}^{\mathbf{q}}\leq c_{10}v^{p-\mathbf{q}}.
\end{eqnarray}
Remind that $\mH(v)\neq\emptyset$,  $\mH(v,z)\neq\emptyset$ whatever $v>0$ and $z\geq 2$, see Remark \ref{rem:existence-of-Lambdas}.  Hence, in view of (\ref{eq251:proof-add}) for any $f$
$$
\Lambda_p(v,f,\mathbf{u})\leq \inf_{z\geq 2}\big[c_{10}v^{p-\mathbf{q}}+z^{-\mathbf{u}}\big]=c_{10}v^{p-\mathbf{q}}.
$$
It remains to note that the right
  hand side of the obtained inequality is independent of $f$ and the second assertion of the theorem follows from this inequality, (\ref{eq222:proof-prop:measure-of-bias-deconv}) and Theorem
\ref{th:consec-th1-unbounded}.

\smallskip

$\mathbf{4^0.}\;$ Since $C_{\mathbf{\infty}}=1$ we obtain in view of  (\ref{eq9:proof-add}) for all $f\in\bB_{\infty,d}(D)$
\begin{eqnarray*}
\label{eq1111:proof-add}
&&\big\|\cB_{\vec{h}}(\cdot,f)\big\|_\infty\leq \big[3(2c_\cK)^d\|\cK\|^d_\infty +1\big]D,\quad\forall \vec{h}\in (0,\infty)^d.
\end{eqnarray*}
It yields for any $\vec{s}\in (1,\infty)$ in view of (\ref{eq2:proof-prop:measure-of-bias-deconv})  if $\mathbf{q}=\infty$
\begin{equation*}
\int_{\cA\big(\vec{h},f,v\big)} \big|\cB_{\vec{h}}(x,f)\big|^{p}\nu_d(\rd x)\leq c_{11} \nu_d\left\{\cA\big(\vec{h},f,v\big)\right\}\leq c_{12}\sum_{j\in\bar{J}(\vec{h},v)} v^{-s_j}\Big[\mathbf{B}_{j,\mathbf{s},\bF}\big(h_{j}\big)\Big]^{s_j}.
\end{equation*}
Since the left hand side of the obtained inequality is independent of $f$ and the left hand side is independent of $\vec{s}$ we conclude that
\begin{eqnarray}
\label{eq122-new:proof-add}
&&\sup_{f\in\bF}\Lambda_{p}(v,f,\mathbf{u})\leq c_{12}\blL_{\vec{s}}(v,\bF,\mathbf{u}),  \quad \forall v>0,\; \vec{s}\in (1,\infty)^d.
\end{eqnarray}
The third  assertion of the theorem  follows now from (\ref{eq222:proof-prop:measure-of-bias-deconv}), (\ref{eq122-new:proof-add}) and Theorem \ref{th:consec-th1-unbounded}.

\smallskip

$\mathbf{5^0.}\;$ We have already seen (Corollary 1), that $B^*_{\vec{h}}(\cdot,f)\leq 2\sup_{\eta\in\cH: \eta\leq h}B_{\vec{\eta}}(\cdot,f)$ if $\vec{h}=(h,\ldots,h)\in\cH^d_{\text{isotr}}$. Therefore by definition of $\cB_{\vec{h}}(\cdot,f)$:
\begin{eqnarray}
\label{eq116:proof-add}
\cB_{\vec{h}}(\cdot,f)&\leq& 5\sup_{\eta\in\cH: \eta\leq h}B_{\vec{\eta}}(\cdot,f)\leq 5\sup_{\eta\in\cH:\eta\leq h}\sum_{j=1}^d\sup_{J\in\mJ}\left[\big|K_{\vec{\eta}}\big|\circ b^*_{\eta,f,j} \right]_{J}(x).
\end{eqnarray}
where, remind $\vec{\eta}=(\eta,\ldots,\eta)\in\cH^d_{\text{isotr}}$.
We remark that (\ref{eq116:proof-add}) is similar to (\ref{eq2001:proof-prop:measure-of-bias-deconv}) but the maximal operator is not involved in this bound. This, in its turn, allows to consider $\vec{s}\in [1,\infty)^d$.

Indeed, similarly to (\ref{eq2000:proof-prop:measure-of-bias-deconv}) we have for any $v>0$,  applying (\ref{eq:Young-restricted}) with $\mathbf{t}=\infty$
\begin{eqnarray}
\label{eq117:proof-add}
\cB_{\vec{h}}(x,f)\leq5\sup_{\eta\in\cH:\eta\leq h}\sum_{j\in\bar{J}(\vec{h},v)}\sup_{J\in\mJ}\left[\big|K_{\vec{\eta}}\big|\circ b^*_{\eta,f,j} \right]_{J}(x)+4^{-1}v, \quad \forall f\in\bF.
\end{eqnarray}
We obtain for any $f\in\bF$, $v>0$ and  $\vec{s}=(s_1,\ldots,s_d)\in [1,\infty)^d$ applying consecutively  the Markov inequality and    (\ref{eq:Young-restricted}) with $\mathbf{t}=s_j$
\begin{eqnarray*}
\label{eq118:proof-add}
\nu_d\left(\cA\big(\vec{h},f,v\big)\right)
\leq c_{13}\sum_{j\in\bar{J}(\vec{h},v)} v^{-s_j}\sum_{\eta\in\cH: \eta\leq h}\big\|b^*_{\eta,f,j}\big\|_{s_j}^{s_j}\leq c_{14}\sum_{j\in\bar{J}(\blh,\bly)} v^{-s_j}\Big[\mathbf{B}^*_{j,s_j,\bF}\big(h\big)\Big]^{s_j}.
\end{eqnarray*}
We note that the obtained inequality coincides with (\ref{eq2:proof-prop:measure-of-bias-deconv}) if one replaces $\mathbf{B}_{j,s_j,\bF}(\cdot)$ by $\mathbf{B}^*_{j,s_j,\bF}(\cdot)$. It remains to remark  that $\mathbf{B}_{j,s_j,\bF}(\cdot)\le \mathbf{B}^*_{j,s_j,\bF}(\cdot)$. Indeed,
$$
b_{\mathbf{v},f,j}(x)=\lim_{k\to\infty}\sup_{h\in\cH:\: e^{-k}\leq h\leq \mathbf{v}}b^*_{h,f,j}(x).
$$
Therefore, by the monotone convergence theorem  and the triangle inequality for any $s\in [1,\infty)$
\begin{eqnarray*}
\mathbf{B}_{j,s,\bF}(\mathbf{h})&:=&\sup_{f\in\bF}\|b_{\mathbf{v},f,j}\|_s=\sup_{f\in\bF}\lim_{k\to\infty}\Big\|\sup_{h\in\cH:\: e^{-k}\leq h\leq \mathbf{h}}b^*_{h,f,j}\Big\|_s
\\
&\leq&\sup_{f\in\bF}\lim_{k\to\infty}\sum_{h\in\cH:\; e^{-k}\leq h\leq \mathbf{h}}\Big\|b^*_{h,f,j}\Big\|_s
=\sup_{f\in\bF}\sum_{h\in\cH:\; h\leq \mathbf{h}}\Big\|b^*_{h,f,j}\Big\|_s=:\mathbf{B}^*_{j,s,\bF}(\mathbf{v}).
\end{eqnarray*}
The fourth  statement of the theorem  follows now from (\ref{eq222:proof-prop:measure-of-bias-deconv}), (\ref{eq22222:proof-prop:measure-of-bias-deconv}), (\ref{eq251:proof-add})
 and Theorem \ref{th:consec-th1-unbounded}.
\epr

\subsection{Proof of Assertion \ref{assert:example-holder}}
\label{sec:subcec:proof-of-assertion}

Obviously $\bF_K\big(\vec{\beta},\vec{L}\big)\subset\bB_{\infty,d}(L_\infty)$. Thus, we can choose $D=L_\infty$ and $\mathbf{q}=\infty$, which implies $\mathbf{u}=\infty$.
For any $v>0$ let $\vec{\blh}(v)=\big(\blh_1(v),\ldots, \blh_d(v)\big)$, where
$$
\blh_j(v)=\max\big\{h\in\cH:\; h\leq \big(\BL L_0 L^{-1}_j v\big)^{1/\beta_j}\big\},\; j=1,\ldots,d,
$$
and $\BL\in (0,1)$ is chosen to satisfy $\BL L_0\leq \mathbf{c}$. This in its turn implies $\BL L_0<1$.

This choice of $\vec{\blh}(v)$ together with the definition of the class $\bF_{K}\big(\vec{\beta},\vec{L}\big)$ implies that
\begin{gather}
\label{eq0:holderian-classes}
J\big(\vec{\blh}(v),v\big)=\{1,\ldots,d\},\quad \forall v>0;
\\
\label{eq1:holderian-classes}
\vec{\blh}(v)\in(0,1]^d,\;\;\;\; \forall v\in(0,\BL^{-1}].
\end{gather}
Moreover, there exists $T_1:=T_1\big(\vec{\beta}\big)<\infty$ {\it independent} of $\vec{L}$ such that
\begin{equation}
\label{eq-new:sum-of-log}
\limsup_{n\to\infty}\;(\ln{n})^{-1}\sup_{v\in V_n}\sum_{j=1}^d\big|\ln{\big(\blh_j(v)\big)}\big|\leq T_1, \quad V_n=\Big[\delta_n^{\frac{\beta(\alpha)}{1+\beta(\alpha)}}, 1\Big].
\end{equation}
Then set $T_2=e^{\frac{d}{2}+\sum_{j=1}^m\mu_j(\alpha)}\sqrt{T_1+2}(\BL L_0)^{-\frac{1}{2\beta(\alpha)}}$.

\vskip0.1cm

We have in view of (\ref{eq1:holderian-classes}) and (\ref{eq-new:sum-of-log})
for all $n$ large enough and any $v\in V_n$
\begin{eqnarray}
\label{eq2:example}
&&F_n\big(\vec{\blh}(v)\big)\leq \frac{\sqrt{(T_1+2)\ln{n}}}{\sqrt{n}\prod_{j=1}^d\big(\blh_j(v)\big)^{\frac{1}{2}+\blg_j(\alpha)}}\leq T_2\sqrt{\delta_n}v^{-\frac{1}{2\beta(\alpha)}}
\\*[1mm]
\label{eq3:example}
&&G_n\big(\vec{\blh}(v)\big)\leq\frac{(T_1+2)\ln{n}}{n\prod_{j=1}^d\big(\blh_j(v)\big)^{1+\blg_j(\alpha)}}\leq \frac{(T_1+2)\ln{n}}{n\prod_{j=1}^d\big(\blh_j(v)\big)^{1+2\blg_j(\alpha)}}\leq T^2_2\delta_nv^{-\frac{1}{\beta(\alpha)}}.
\end{eqnarray}
 Setting $T_3=(\sqrt{2}a^{-1}T_2)^{\frac{2\beta(\alpha)}{2\beta(\alpha)+1}}$ and $T_4=(a^{-1}T^2_2)^{\frac{\beta(\alpha)}{\beta(\alpha)+1}}$ we obtain in view of (\ref{eq2:example}) and (\ref{eq3:example}) for all $n$ large enough
\begin{eqnarray}
\label{eq4:example}
&& \vec{\blh}(v)\in\mH(v),\;\;\forall v\in\Big[T_4\delta_n^{\frac{\beta(\alpha)}{1+\beta(\alpha)}}, 1\Big],\,\quad \vec{\blh}(v)\in \mH(v,2),\;\;\forall v\in\big[T_3\varphi_n, 1\big].
\end{eqnarray}
It is worth noting that $T_2>1$, which implies $T_4>1$, and $\delta_n^{\frac{\beta(\alpha)}{1+\beta(\alpha)}}\varphi^{-1}_n\to 0, n\to\infty$. Choose
$$
\underline{\Bv}=T_4\delta_n^{\frac{\beta(\alpha)}{1+\beta(\alpha)}},
\quad \overline{\Bv}=T_3\varphi_n.
$$
Since $\vec{\blh}(v)\in\mH(v)$ for any $v\in[\underline{\Bv},\overline{\Bv}]$ in view of  (\ref{eq4:example}), we deduce from (\ref{eq0:holderian-classes})
and (\ref{eq2:example}) for any $\vec{s}$
$$
\blL_{\vec{s}}\big(v,\bF_K\big(\vec{\beta},\vec{L}\big)\big)\leq v^{-2}F_n^2\big(\vec{h}(v)\big)\leq T_2^2\delta_nv^{-2-1/\beta(\alpha)} ,\quad\forall v\in[\underline{\Bv},\overline{\Bv}].
$$
This, in its turn, yields for any $\vec{s}$
\begin{eqnarray}
\label{eq4000:example}
&&\int_{\underline{\Bv}}^{\overline{\Bv}} v^{p-1}\blL_{\vec{s}}\big(v,\bF_K\big(\vec{\beta},\vec{L}\big)\big)\rd v\leq T_5\delta_nZ(\underline{\Bv},\overline{\Bv}),
\end{eqnarray}
where we have denoted $T_5=T_2^2\big\{1\vee |p-2-1/\beta(\alpha)|^{-1}\big\}$ and
$$
Z(\underline{\Bv},\overline{\Bv})=\overline{\Bv}^{p-2-1/\beta(\alpha)}\mathrm{1}_{\{p>2+1/\beta(\alpha)\}}+\underline{\Bv}^{p-2-1/\beta(\alpha)}\mathrm{1}_{\{p<2+1/\beta(\alpha)\}}+
\ln{(\overline{\Bv}/\underline{\Bv})}\mathrm{1}_{\{p=2+1/\beta(\alpha)\}}.
$$
Moreover, since $\vec{\blh}(\overline{\Bv})\in\mH(\overline{\Bv},2)$ in view of (\ref{eq4:example}), we deduce from (\ref{eq0:holderian-classes}) that
 for any $\vec{s}$
\begin{eqnarray}
\label{eq5:example}
&&\blL_{\vec{s}}\big(v,\bF_K\big(\vec{\beta},\vec{L}\big),\infty\big)=0.
\end{eqnarray}
At  last, putting $T_6=T_4^{p-1}(1+\ln{T_4})^{t(\bH)}$, we  obtain
\begin{eqnarray}
\label{eq6:example}
l_\bH(\underline{\Bv}):=\underline{\Bv}^{p-1}(1+|\ln{(\underline{\Bv})}|)^{t(\bH)}\leq T_6\delta_n^{\frac{(p-1)\beta(\alpha)}{1+\beta(\alpha)}}(1+\ln{n})^{t(\bH)}.
\end{eqnarray}
Applying the third assertion of Theorem \ref{th:minimax-abstract}, we deduce from (\ref{eq4000:example}), (\ref{eq5:example}) and (\ref{eq6:example}) that
\begin{equation*}
\cR^{(p)}_n[\widehat{f}_{\vec{\mathbf{h}}(\cdot)}, f]\leq C^{(1)}\bigg[T_6\delta_n^{\frac{(p-1)\beta(\alpha)}{1+\beta(\alpha)}}(1+\ln{n})^{t(\bH)}+
T_5\delta_nZ(\underline{\Bv},\overline{\Bv})
\bigg]^{\frac{1}{p}}+\mathbf{C}_pn^{-\frac{1}{2}}.
\end{equation*}
After elementary  computations we come to the statement of Assertion \ref{assert:example-holder}.
\epr

\bibliographystyle{agsm}

\end{document}